\documentclass[final,english]{article}

\usepackage[OT1]{fontenc}
\usepackage[latin1]{inputenc}
\usepackage{graphicx}
\usepackage{amsmath,amssymb}
\usepackage[normalem]{ulem}

\makeatletter
\newtheorem{thm}{Theorem}
\newtheorem{prop}[thm]{Proposition}

\newtheorem{defin}[thm]{Definition}
\newtheorem{lem}[thm]{Lemma}
\newcommand{\proof}{}

\usepackage{epsfig}
\usepackage{latexsym}
\usepackage{tabularx}
\usepackage{fancyhdr}
\usepackage{pifont}    
\usepackage{babel}
\usepackage{tikz}
\usepackage{pgfplots,pgfplotstable}

\usepackage{color}

\bigskip

\def\RR{\mathbb{R}}
\def\ZZ{{\mathbb{Z}}}
\def\NN{\mathbb{N}^{*}}

\def\AA{\mathcal{A}}

\def\HH{{\mathcal{H}}}
\def\AA{{\mathcal{A}}}

\def\XX{{\mathcal{X}}}

\def\vv{{\bf v}}

\newcommand{\FG}[1]{#1}
\newcommand{\integer}[2]{[\hspace{-0.40ex} [ {#1},{#2} ]\hspace{-0.35ex} ]}

\renewcommand{\div} {\mathop{\rm div}\nolimits}

\makeatother  

\title{The case of Neumann, Robin and periodic lateral condition for the semi infinite generalized Graetz problem and applications.}

\author{
Valention Debarnot, J\'er\^ome Fehrenbach, Fr\'ed\'eric de Gournay, L\'eo Martire
\thanks{
Institut de Mathématiques de Toulouse (UMR 5219), Université de Toulouse, CNRS 
UPS, F-31062 Toulouse, France F-31077 Toulouse, France 
INSA, F-31077 Toulouse, France
({\tt frederic@degournay.fr}), ({\tt jerome.fehrenbach@math.univ-toulouse.fr})
 } 
 }
\begin{document}

\maketitle

\begin{abstract}
The Graetz problem is a convection-diffusion equation in a pipe invariant along a direction. The contribution of the present work is to propose a mathematical analysis of the Neumann, Robin and periodic boundary condition on the boundary of a semi-infinite pipe. The solution in the 3D space of the original problem is reduced to eigenproblems in the 2D section of the pipe. The set of solutions is described, its structure depends on the type of boundary condition and of the sign of the total flow of the fluid. This analysis is  the cornerstone of numerical methods to solve Graetz problem in finite pipes, semi infinite pipes and exchangers of arbitrary cross section. Numerical test-cases illustrate the capabilities of these methods to provide solutions in various configurations.
\end{abstract}

\section{Introduction}
\subsection{Context}
The seminal work of Graetz in the late 19th century adressed a stationnary  convection-diffusion problem inside an axi-symmetrical cylindrical pipe \cite{graetz1885warmeleitungsfahigkeit}, where the regime was supposed to be convection-dominated which means that the longitudinal diffusion \FG{was} neglected. It was the first contribution to the modelling of convective transport coupled with diffusion, with important applications nowadays as the parallel convective exchangers involved in  heating or cooling systems \cite{shah_Dusan_Sekulic_book}, haemodialysis \cite{gostoli_80}, 
and heat exchangers  \cite{kragh_07}.  
The first extension to the Graetz problem, known as the ``extended Graetz problem'' takes into account longitudinal diffusion \cite{michelsen1974graetz,ebadian1989exact,weigand2001extended,LAHJOMRI_02}.
Papoutsakis {\em et al.} in \cite{papoutsakis_80a,papoutsakis_80b} introduced a symmetric operator acting on a two-components space that solves the  extended Graetz problem in axi-symmetrical configurations. The  so-called ``conjugated Graetz problem'' where multiple solid or  fluid phases are taken into account was proposed in 
\cite{papoutsakis_81_1,papoutsakis_81_2} in the case of an axi-symmetrical configuration.
These successive models aimed at taking into account more and more complex and realistic situations, and when only axi-symmetrical configurations were considered the equations boiled down to one-dimensional problems. 
The adaptation to \FG{parallel} plates heat exchangers  of these one-dimensional models, together with a parametric study  was proposed in \cite{ho1998analytical}. The reader may also consult \cite{dorfman2009conjugate} for a review on the conjugated Graetz problem.

The work on non-axisymmetrical configurations  was initiated in  \cite{pierre_C_09} where the operator was proved to be self-adjoint  with compact resolvent when Dirichlet boundary conditions are applied on the boundary of the domain. In the case of a single fluid stream the negative eigenvalues correspond to downstream propagation, and positive eigenvalues to upstream propagation. The main novelty was that arbitrary geometries were adressed, and a detailed mathematical analysis of the Dirichlet problem was proposed.
Numerical methods for the approximation of this operator and error estimates where provided in \cite{fehrenbach2012generalized}.

The objective of the present work is to extend the work of \cite{fehrenbach2012generalized} and  provide explicit methods with general  lateral boundary conditions, beyond the Dirichlet case. The cross section of the domain has an arbitrary geometry and can incorporate different fluid domains, possibly with opposite signs of the velocity.
The lateral boundary conditions that we address can be Dirichlet, Neumann, Robin, periodic or a mixture of these different cases on different parts of the boundary.
The periodic boundary conditions with rectangular or hexagonal cell are adapted to the analysis of micro-exchangers, where a design pattern is repeated.

\subsection{Setting}
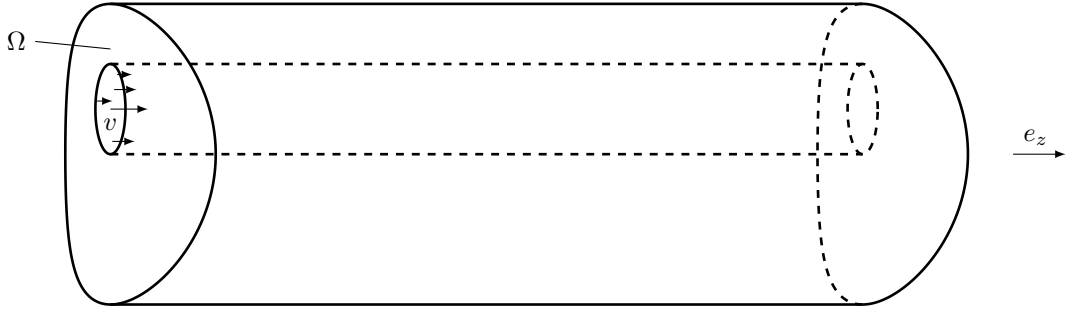
\begin{figure}
\begin{center}
\begin{tikzpicture}
\draw  [domain=0:360,scale=2,variable=\t,samples=200,line width=1pt] plot({0.5*cos(\t)+0.2*cos(\t)*cos(\t)}, {sin(\t)});
\draw  [domain=90:270,scale=2,variable=\t,samples=200,line width=1pt,dashed] plot({5+0.5*cos(\t)+0.2*cos(\t)*cos(\t)}, {sin(\t)});
\draw  [domain=270:450,scale=2,variable=\t,samples=200,line width=1pt] plot({5+0.5*cos(\t)+0.2*cos(\t)*cos(\t)}, {sin(\t)});
\draw [line width=1pt] (0,2)--(10,2);
\draw [line width=1pt] (0,-2)--(10,-2);
\draw  [domain=0:360,scale=2,variable=\t,samples=200,line width=1pt] plot({0.1*cos(\t)}, {0.3+0.3*sin(\t)});
\draw  [domain=0:360,scale=2,variable=\t,samples=200,line width=1pt,dashed] plot({5+0.1*cos(\t)}, {0.3+0.3*sin(\t)});
\draw [line width=1pt,dashed] (0,1.2)--(10,1.2);
\draw [line width=1pt,dashed] (0,0)--(10,0);
\draw [>=latex,->] (0.05,0.86)--(0.35,0.86);
\draw [>=latex,->] (0.025,0.17)--(0.325,0.17);
\draw [>=latex,->] (0,0.6)--(0.5,0.6);
\draw [>=latex,->] (0.09,1.06)--(0.29,1.06);
\draw [>=latex,->] (-0.1764,0.7073)--(0.03,0.7073);
\draw [>=latex,->] (12,0)--(12.7,0);
\draw (0,0.6) node [below] {$v$}; 
\draw (12.3,0) node [above] {$e_z$}; 
\draw (-1,1.5) node [left] {$\Omega$}; 
\draw (-1,1.5)--(0,1.4);
\end{tikzpicture}
\end{center}
\caption{The domain $\Omega\times I$ where the Graetz problem is posed.}\label{fig;omegaI}
\end{figure}

In convection-dominated heat or mass transfer, we address the generalized Graetz problem which occurs in a cylinder of arbitrary section $\Omega$ and of length $I$, possibly $I=\RR^+$, see Figure \ref{fig;omegaI}. The diffusion coefficient is supposed to be invariant by translation along $e_z$ the axis of the cylinder. Similarly, the velocity vector $v$ is supposed to be oriented in the direction of the axis of the cylinder, that is $v=he_z$ with $h\in L^\infty(\Omega)$. The equation for the temperature $T$ inside the domain is then
\begin{equation}
\label{eq::inside}
\tag{E} c\partial_{zz} T +\div(\sigma \nabla T) -h\partial_z T=0 \text{ on } \Omega \times I, 
\end{equation}
with diffusion coefficients $c,\sigma >0$ bounded in $\Omega$ with bounded inverse.
The lateral homogeneous boundary conditions~\eqref{eq:lateralboundary}  may be of Neumann, Dirichlet, Robin and periodic type, respectively on $\Gamma_N$,$\Gamma_D$,$\Gamma_R$,$\Gamma_\sharp \subset \partial \Omega$  given by 
\begin{equation}
\tag{LBC}
\begin{cases} 
\sigma \nabla T\cdot n=0 \text{ on } \Gamma_N \times I \text{: Neumann, and/or}\\
T=0 \text{ on } \Gamma_D \times I \text{: Dirichlet, and/or}\\
\sigma \nabla T \cdot n+aT=0 \text{ on } \Gamma_R \times I \text{: Robin, and/or} \\
T  \text{ is periodic on } \Gamma_\sharp \times I \text{: periodic} \\
\end{cases} \label{eq:lateralboundary}
\end{equation}
where $a> 0$ in the Robin condition, and $\Gamma_\sharp$ must be taylored to support periodic conditions (e.g. $\Omega$ is the unit square, $\Gamma_\sharp=(\{x=0\}\cup\{x=1\})\cap \partial \Omega$ and the boundary condition is $T(0,y)=T(1,y)$). As usual, the \FG{$\Gamma$'s} involved in the definition of the boundary condition must form a partition of $\partial \Omega$. Note that the Neumann (resp. Dirichlet) boundary conditions are degenerate cases of the Robin condition corresponding to $a=0$ (resp. $a=+\infty$). The
Inlet/Outlet boundary condition~\eqref{eq:ioboundary} is  of  Dirichlet and/or of Neumann type and is given by 
\begin{equation}
\tag{I/OBC}
T=T_D \text{ on } \Omega_D\text{ and }\partial_zT=S_N \text{ on } \Omega_N \text{ with } \Omega_D\cup \Omega_N=\Omega\times \partial I.
\label{eq:ioboundary}
\end{equation}
In the case $I=\RR^+$, we intentionnally stay vague about the definition of $\partial I$, it is one of the results of this \FG{work} to determine whether an (I/OBC) is needed on $z=+\infty$.

A more realistic model in regimes of high velocities takes into account a viscosity term, see e.g.~\cite{jeong2006extended} where a study in a microchannel including viscous effects and longitudinal conduction is performed. 
Our approach can also account for viscosity, the details are presented in Section~\ref{ssec:viscous}.

\subsection{Lax-Milgram}
\label{sec:lax-milgram}
Note that the equation \eqref{eq::inside} is an elliptic equation with an additionnal convective term. It is possible to use
 Lax-Milgram's theorem \cite{rudin1991functional} under the hypothesis that the Inlet/Outlet boundary condition is  Dirichlet in the region where the flow is incoming. More precisely:
\begin{prop}
Let $I=[z_1,z_2]$ and $\omega_\pm=\{x\text{ s.t } \pm h(x)> 0\}$. If 
\[ \omega_+\times\{z_1\}\subset\Omega_D \text{ and }\omega_-\times\{z_2\}\subset\Omega_D\]
and if $T_D$ and $S_N$ are regular enough,
then there exists a unique solution to \eqref{eq::inside} with the boundary conditions \eqref{eq:lateralboundary} and \eqref{eq:ioboundary}.
\end{prop}

The proof is only sketched here for the sake of completness. \FG{Denote $\XX$ the} natural space of elements where the solution is sought, that is
\[\XX=\{T \in H^1(\Omega\times I) \text{ s.t. } T=0 \text{ on } \left(\Gamma_D\times I\right)\cup \Omega_D \text{ and } T \text{ periodic on }\Gamma_\sharp\times I\}.\]
Non-homogeneous Dirichlet boundary conditions of \eqref{eq:ioboundary} are solved using a lift of $T_D$, still denoted  $T_D$   that satisfies the lateral boundary conditions \eqref{eq:lateralboundary} with $\partial_z T_D =0$ on $\Omega_N$ and denote 
\[f_D=c\partial_{zz}T_D +\div(\sigma \nabla T_D)-h\partial_zT_D.\]
The change of unknown $\widetilde T=T-T_D$ where $T$ solves \eqref{eq::inside} and \eqref{eq:lateralboundary}, leads to the following variational formulation: find $\widetilde T\in \XX$ such that for every $\phi\in \XX$: 
$$\underbrace{\int_{\Omega\times I}c\partial_{z}\widetilde T \partial_z\phi+\sigma \nabla \widetilde T \cdot \nabla\phi+h\partial_z\widetilde T \phi + \int_{\Gamma_R} a \widetilde T\phi}_{{ b(\widetilde T,\phi)}}+\underbrace{\int_{\Omega_N}S_N \phi-\int_{\Omega\times I} f_D\phi}_{{\ell(\phi)}}=0.$$

The term $b(T,\phi)$ is bilinear in $(T,\phi)$ and continuous for the standard norm of $\XX$, the term $\ell(\phi)$ is linear continuous if $T_D$ and $S_N$ are regular enough. It remains to study the coercivity of $b$. 
\[b(T,T)=\int_{\Omega\times I}c\partial_{z}T.\partial_zT+\sigma\nabla T\cdot\nabla T+h\partial_zT.T=\int_{\Omega\times I}\|\nabla_{\rm 3D}T\|_\kappa^2+\dfrac{1}{2}\int_{\Omega\times I}h\partial_z(T^2),
\]
where $\kappa$ is a positive matrix with diagonal entries $(\sigma,\sigma,c)$ in the basis $(e_x,e_y,e_z)$.
The first term is coercive. The second term is 
$$\dfrac{1}{2}\int_\Omega hT^2|_{z=z_1}^{z=z_2}=\dfrac{1}{2}\int_{\Omega\times\{z_2\}} hT^2-\dfrac{1}{2}\int_{\Omega\times\{z_1\}} hT^2.$$
It is nonnegative for all $T\in\XX$ if and only if the Inlet/Outlet condition is of Dirichlet type at the boundary where the flow is entering the domain ($z=z_1$ if $h>0$, and $z=z_2$ if $h<0$). 

\subsection{Presentation of the paper}
The objective of the present paper is to provide a general framework that allows to solve \eqref{eq::inside} with any type of boundary condition beyond the case where Lax-Milgram's theorem can be used.
Section~\ref{sec:position} details the notation and the main properties of the operator involved in the solution, as well as the modifications required to take into account a viscosity term.
The main results, namely Theorems \ref{theo:invertibility} and~\ref{th:decomposition}, are detailed in Section~\ref{sec:decomp}, their proof is postponed to the Appendix. In Section~\ref{sec:semi-infinite} we solve the problem in a semi-inifinite domain and show that depending on the case the temperature at infinity $T_\infty$ can either be a free parameter of the problem or be imposed by the other condition. 
In Section~\ref{sec:finite} we adress the case of a domain of finite length, and numerical strategies are detailed in the different cases depending on the lateral boundary condition and on the Inlet/Outlet condition. Test cases are presented in Section~\ref{sec:numerical}.

\section{State of the art and position of the problem}

\label{sec:position}
The equation \eqref{eq::inside} may be interpreted as an evolution equation in the variable $z$ if it is cast into
\begin{equation}
\label{eq:def:A}
\partial_z \begin{pmatrix} \partial_z T \\ T\end{pmatrix}={\mathcal A} \begin{pmatrix} \partial_z T \\ T\end{pmatrix} \text{ on } \Omega\times I\text{, with } {\mathcal A}\begin{pmatrix} u \\ s\end{pmatrix}=\begin{pmatrix} hc^{-1}u -c^{-1}\div\sigma\nabla s \\ u\end{pmatrix}. 
\end{equation}
The goal of this section is to guide the reader to the  analysis of \eqref{eq:def:A} that was proposed in  \cite{fehrenbach2012generalized}, to enlarge the frame to Neumann and periodic lateral boundary condition, and to define the notation and state the results that will be used in the sequel. Since $\mathcal A$ is  a symmetric operator with a compact resolvent, classical eigendecomposition leads to \FG{ an explicit representation of the solution of \eqref{eq:def:A} in the basis of eigenvectors} (see e.g. \cite{rudin1991functional}).

\begin{defin}
We say that {\em ``the constants are not controlled''} when $\Gamma_D\cup\Gamma_R=\emptyset$, in other words when there is no Dirichlet or Robin condition on the lateral part of the boundary of the domain. The case where the constants are not controlled and in addition $\int_\Omega h=0$ is called the {\em ``balanced case''}. 
\end{defin}

As we prove in this section, the case where the constants are not controlled is a case where the constants are a solution of \eqref{eq::inside} and the balanced case is a case where $\mathcal A$ admits a non-trivial kernel.
\subsection{Study of the operator $\mathcal A$}
\label{sec:operatorA}
In this section we detail the Hilbert space, the scalar product, the kernel, range and pseudo-inverse of the symmetric operator $\mathcal A$. 

{\bf Hilbert space and scalar product.} First, introduce the space $H$ that encodes the lateral boundary condition.
 When the constants are controlled define:
\[H=\{s \in H^1(\Omega), \text{ such that } s=0 \text{ on } \Gamma_D \text{ and } s \text{ periodic on } \Gamma_\sharp\}.\]
If there is no Dirichlet or Robin boundary condition, hence no control on the constants, quotient by the constants and define: 
\[H=\{s \in H^1(\Omega)/\RR, \text{ such that }  s \text{ periodic on } \Gamma_\sharp\}.\]
Then, define the Hilbert space $\HH$ as 
\[\HH=\{(u,s)\,|\,u\in L^2(\Omega),s\in H\}\]
which is endowed with the scalar product:
$$\left(\begin{pmatrix} u, s\end{pmatrix}|\begin{pmatrix} u',  s'\end{pmatrix}\right)_\HH=\int_\Omega cuu'+\sigma \nabla s\cdot \nabla s'+  \int_{\Gamma_R} a ss'.$$
The crucial step in showing that $\HH$ is a Hilbert space is to show that the scalar product is definite. Setting  $\left(\begin{pmatrix} u, s\end{pmatrix}|\begin{pmatrix} u,  s\end{pmatrix}\right)_\HH=0$ immediatly gives $u=0$ and $\nabla s=0$, hence $s$ is a constant. If the constants are controlled, then $\Gamma_D\cup\Gamma_R\ne\emptyset$ and $s=0$, whereas if the constants are not controlled then $s$ is a constant and $s=0$ in $H$.

The domain of the operator $\mathcal A $ is:
\[\mathcal{D(A)}=\{(u,s)\in \HH, u \in H^1(\Omega), \div(\sigma \nabla s) \in L^2(\Omega) +\text{ boundary conditions \eqref{eq:lateralboundary}}\},\]
where the boundary conditions are $u\in H$, and $\sigma\nabla s\cdot n$ is equal to $0$ on $\Gamma_N$, is equal to $-as$ on $\Gamma_R$ and is periodic on $\Gamma_\sharp$.
On $\mathcal{D(A)}$, the operator is symmetric as we prove now. Let $\phi=(u,s)$ and  $\phi'=(u',s') \in \mathcal{D(\mathcal A)}$:
\begin{eqnarray*}
&&\left({\mathcal A}\phi|\phi'\right)_\HH=\int_\Omega (hu-\div \sigma \nabla s)u'+\sigma \nabla u\cdot \nabla s'+\int_{\Gamma_R} a us'
\\
&=&\int_\Omega huu'+\sigma \nabla u\cdot \nabla s'+\sigma \nabla u'\cdot \nabla s+\underbrace{\int_{\partial \Omega}(-\sigma \nabla s \cdot n)u' +\int_{\Gamma_R} a us'}_{\bf (1)},
\end{eqnarray*}
and the term $(\bf 1)$ is symmetric thanks to \eqref{eq:lateralboundary} on $\mathcal{D(A)}$.

{\bf Inverse of the Laplacian}
Define the inverse of the Laplace operator as:
\[u=\Delta_\sigma^{-1}f \text{ iff } \left\{
\begin{array}{l}
\div(\sigma \nabla u)=f, \text{ and } \\
u\in H \\
+ \text{ boundary conditions},\end{array}
\right.
\] 
where the boundary conditions are $\sigma \nabla u \cdot n=0$ on $\Gamma_N$ and $\sigma \nabla u \cdot n+au=0$ on $\Gamma_R$.
If the constants are controlled, then $\Delta_\sigma^{-1}$ is well defined on $L^2(\Omega)$, whereas if there is only Neumann or periodic boundary conditions  (no control of the constants), the operator $\Delta_\sigma^{-1}$ is only defined if $f\in L^2_m(\Omega)$, the subspace of $L^2(\Omega)$ with null average. 

{\bf Kernel of $\mathcal A$}
Following from the definition of $\mathcal A$ in \eqref{eq:def:A}, the kernel of $\mathcal A$ is the set of $(u, s)$ in $\mathcal D(\mathcal A)$ such that
$$ u= 0 \text{ in } H \text{ and } hu-\div(\sigma \nabla s)=0.$$
When the constants are controlled, both $u$ and $s$ are then equal to $0$. When the constants are not controlled, since $u$ is a constant, then $s=u\Delta_\sigma^{-1}h$ in $\Omega$ which admits a solution \FG{if and only if} $\int_\Omega h=0$.
To summarize the kernel of $\mathcal A$ is:
\[ \mathcal K(\mathcal A)=\begin{cases} Vect(\phi_0=(1,\Delta_\sigma^{-1}h)) \text{ in the balanced case },\\
\{0\} \text{ in the other cases}
\end{cases}\]

{\bf Range and inverse of $\mathcal A$} The range of $\mathcal A$, denoted  $\mathcal R( \mathcal A)$, is defined as the orthogonal of  $\mathcal K(\mathcal A)$ in $\mathcal H$ and the inverse of $\mathcal A$ is an operator from $\mathcal R(\mathcal A)$ to $\mathcal D(\mathcal A)$, defined as follows:
\[\forall \phi=(u,s) \in \mathcal R( \mathcal A), \mathcal A ^{-1}\phi=\begin{cases} (s,\Delta_\sigma^{-1}(hs-cu)) \text{ if the constants are controlled} \\
(s+k,\Delta_\sigma^{-1}(hs-cu+hk)), k\in \RR \text{ if not}.
\end{cases}.\]
When the constants are not controlled, the constant $k\in \RR$ is chosen so that
\[\begin{cases}\int_\Omega hs-cu+hk=0 \text{ in the non-balanced case}, \\
(\mathcal A^{-1}\phi,\phi_0)_{\mathcal H}=0 \text{ in the balanced case}.
\end{cases}\]
It is easily checked that for all $\phi \in \mathcal R(\mathcal A)$, $\mathcal A^{-1} \phi \in \mathcal D(\mathcal A)$ and that $\mathcal A\mathcal A^{-1}\phi=\phi$. The operator $\mathcal A^{-1}$ is then symmetric (as a consequence of the symmetry of $\mathcal A$).
Note also that in the balanced case, one can also write:
\[\forall \phi=(u,s) \in \mathcal R( \mathcal A),\quad   \mathcal A^{-1}\phi=(s,\Delta_\sigma^{-1}(hs-cu)) + k\phi_0.\]

{\bf Eigenvalue decomposition of $\mathcal A$}
The operator $\mathcal A^{-1}$ is a compact self-adjoint operator on $\mathcal R(\mathcal A)$. To prove this let $\phi_n=(u_n,s_n)$ be a bounded sequence in $\HH$. Then up to a subsequence it is a weakly convergent sequence and $s_n$ converges strongly in $L^2(\Omega)$. Using the fact that $\Delta_\sigma^{-1}$ is a compact operator from $L^2$ to $H$ finishes the proof.  
We denote by $\lambda_i$ the non-zero ordered eigenvalues of ${\mathcal A}$ and by $\phi_i=(U_i,\lambda_i^{-1}U_i)$ the corresponding eigenvectors. By convention, $\lambda_i$ is of the sign of $i$ so that
\[-\infty \leftarrow \lambda_{-n}\le \lambda_{-n-1} \le \dots \le \lambda_{-1}<0<\lambda_1\le \dots \le \lambda_{n-1}\le \lambda_{n} \rightarrow +\infty. \]

In the balanced case, we add to the family $(\phi_i)_i$ the vector $\phi_0=(1,\Delta_\sigma^{-1}h)$, so that the Hilbert space  $\mathcal H$ is the space spanned by the eigenvectors $(\phi_i)_{i\in \ZZ}$. 


\subsection{Solution of the evolution equation}
\label{sec:solevol}
The diagonalization of the operator $\mathcal A$ allows to solve the evolution equation \eqref{eq::inside}:
\begin{equation}
\tag{\ref{eq::inside}}
c\partial_{zz} T +\div(\sigma \nabla T) -h\partial_z T=0 \text{ on } \Omega \times I.
\end{equation}
Let $T \in C^1(I,L^2(\Omega))\cap C^0(I,H)$ be a solution of this equation with corresponding lateral boundary conditions \eqref{eq:lateralboundary}. If we denote $\phi:z\mapsto(\partial_zT(z),T(z))$ in $C^0(\mathcal{H})$, then the equation \eqref{eq::inside} is equivalent to $\partial_z \phi =\mathcal A \phi$, and the solution $\phi$ is given by
\begin{equation}
\label{eq::evolution-spectral}
\phi(z)=\sum_{i\in \ZZ}\frac{(\phi(0)|\phi_i)_{\mathcal H}}{\Vert \phi_i\Vert^2_{\HH}} e^{\lambda_i z}\phi_i.
\end{equation}

One can either identify the first coordinate and integrate w.r.t. $z$ or identify the second coordinate and denote 
\[\psi=\sum_{i\in \ZZ^{*}} \lambda_i^{-1}\frac{ (\phi(0)|\phi_i)_{\mathcal H}}{\Vert \phi_i\Vert_\HH^2} \phi_i,\]
to obtain
\begin{eqnarray*}
T(z)&=&\sum_{i\in \ZZ^{*}}(\psi|\phi_i)_\HH U_i e^{\lambda_i z} \text{ if the constants are controlled, i.e.} \Gamma_D\cup \Gamma_R\ne \emptyset, \\
T(z)&=&\sum_{i\in \ZZ^{*}}(\psi|\phi_i)_\HH U_i e^{\lambda_i z} +a_0\text{ with } a_0 \in \RR
 \text{ if }\Gamma_D\cup \Gamma_R=\emptyset \text{ and } \int_\Omega h \ne 0, \\
T(z)&=&\sum_{i\in \ZZ^{*}}(\psi|\phi_i)_\HH U_i e^{\lambda_i z} +a_0 +a_1(z+\Delta_\sigma^{-1}h) \text{ in the balanced case }\\
&&\text{ with } a_0 \in \RR \text{ and } a_1=\frac{(\phi(0)|\phi_0)_{\mathcal H}}{\Vert \phi_0\Vert^2}.
\end{eqnarray*}

If $\partial_z T$ and $T$ are given at $z=0$ such that $\phi(0)=(\partial_zT(0),T(0))$ belongs to $\mathcal{H}$, then  $\psi$ is uniquely determined. Moreover the constant $a_0$ is also determined by $T_{|z=0}$ (and also $a_1$ in the balanced case). We stress that this solution may not be defined everywhere, indeed the series on the right-hand side of \eqref{eq::evolution-spectral} has to be convergent in some sense and the convergence of the series for $z=0$ is not sufficient to ensure the convergence for $z\ne 0$ due to the multiplication by $e^{\lambda_i z}$ for non-zero $\lambda_i$'s. The set of initial datum $\phi$ that \FG{allows} this series to exist is known as the set of compatible initial condition for the Cauchy problem.

\subsection{Including a viscous term}
\label{ssec:viscous}

Let us consider the following modification of the equation \eqref{eq::inside} where a viscous term is added:
\begin{equation}\label{eq:visco}
c\partial_{zz}T+\div(\sigma\nabla T)-h\partial_z T=\mu|\nabla h|^2.
\end{equation}
\begin{prop}
Let $T$ be the solution of the Graetz equation with viscosity \eqref{eq:visco}. Then there exists an explicit change of unknown function that transforms the problem with viscosity into a problem without viscosity of the form \eqref{eq::inside}. Therefore the solution of the problem \eqref{eq:visco} reduces to the solution of the original problem \eqref{eq::inside}.
\end{prop}
\proof
Once \FG{a} particular solution $\widetilde T$ is found, the change of variable $\widehat T=T-\widetilde T$ transforms by linearity the problem with viscosity \eqref{eq:visco} into the problem without viscosity. We  distinguish different cases, depending on if the constants are controlled or not, and in the case the constants are not controlled we treat separately the non-balanced and the balanced case.  In each case we provide an explicit particular solution $\widetilde T$.

\begin{enumerate}
\item[a)] If the constants are controlled, a particular solution is given by
$$\widetilde T=\Delta_\sigma^{-1}(\mu|\nabla h|^2).$$
\item[b)] If the constants are not controlled, in the non-balanced case 
$$\widetilde T=\alpha z+ \Delta_\sigma^{-1}(\mu|\nabla h|^2+\alpha h),$$
where $\alpha \in \RR$ satisfies
$$\int_\Omega\left(\mu|\nabla h|^2+\alpha h\right)=0.$$
\item[c)] If the constants are not controlled, in the balanced case, the particular solution is given by
$\widetilde T=\alpha (\displaystyle\frac{z^2}{2}+z\Delta_\sigma^{-1} h)+\Delta_\sigma^{-1}\gamma,$ with $\alpha\in \RR$ and $\gamma \in L^{2}(\Omega)$ such that:
\begin{equation}\label{eq:balancedvisco}
\left\{
\begin{array}{l}
\alpha\left(\displaystyle \int_\Omega c-h\Delta_\sigma^{-1} h\right)=\displaystyle  \int_\Omega\mu|\nabla h|^2 \\
\gamma=\mu|\nabla h|^2-\alpha(c-h\Delta_\sigma^{-1} h)\\
\end{array}
\right.
\end{equation} 
The choice of $\alpha$ ensures that $\gamma$ has zero average so that $\Delta_\sigma^{-1}\gamma$ is well defined. Note that $\alpha$ is well defined since
$$\int_\Omega c-h\Delta_\sigma^{-1} h=\int_\Omega c+ \int_\Omega \sigma |\nabla h|^2> 0.$$
Note that the last term is equal to $\Vert \phi_0\Vert^2_{\HH}$.$\Box$
\end{enumerate}

\section{Main decomposition theorem}
\label{sec:decomp}
In this section, the decomposition of a temperature field on the \FG{non-positive eigenspace} is studied. The result stated in Theorem~\ref{th:decomposition} considers different cases depending on the control of constants and the sign of the total flow.
\subsection{Notation and statement of the problem}
The $D(\mathcal A^{\alpha})$ norm or ``$\alpha$-norm'' in short is defined by
\[\Vert \phi\Vert^2_\alpha = \sum_{i\in \ZZ} \lambda_{i}^{2\alpha}\frac{(\phi_i|\phi)_{\mathcal H}^2}{\Vert \phi_i\Vert_{\HH}^2} \quad \forall \phi \in \mathcal H.\]
The space $D(\mathcal A^{\alpha})$ is the set of $\phi \in \mathcal R(\mathcal A)$ whose $\alpha$-norm is $<+\infty$. It is easy to check that $D(\mathcal A^{1})$=$D(\mathcal A)$ and that $D(\mathcal A^{0})$=$\mathcal R(\mathcal A)$.
Define $P$ an orthogonal projection on $\mathcal{H}$ as: 
\begin{equation}
\label{eq:def:P}\forall \phi=(u,s) \in \mathcal{H}, \quad P\phi=(u,0).
\end{equation}
For any $I$ subset of $\ZZ$ define $\pi_I$ the orthogonal projection
\begin{equation}
\label{eq:def:pi}
\pi_I \phi = \sum_{i\in I} \frac{(\phi_i|\phi)_{\mathcal H}}{\Vert \phi_i\Vert_\HH^2} \phi_i
\end{equation}
We denote $\pi_+=\pi_{\NN}$, $\pi_-=\pi_{-\NN}$,  $\pi_0=\pi_{\{0\}}$ and $\mathcal R(\pi_I)=\pi_I(\mathcal{H})$.

The problem of decomposition of a temperature field on the \FG{non-positive eigenspace} is stated as follows:
 
{\em 
For any $\phi \in \mathcal{H}$, find $\psi$ such that 
\begin{equation}
\label{equa::mainpbm2}
P\psi=P\phi \text{ and } \pi_+\psi=0.
\end{equation}}
A similar problem of decomposition on the \FG{non-negative eigenspace} is obtained by replacing $\pi_+$ by $\pi_-$.  All the results of the present section have a counterpart obtained by changing the sign of $z$.

\subsection{Necessary and sufficient condition}
\label{sec:decomp:necessary}
In order to tackle problem \eqref{equa::mainpbm2}, we first consider the following related problem:
\begin{equation}
\label{eq::newproblem}
\text{Find } \psi \in \mathcal R(\pi_-) \text{ such that } \pi_- P\pi_- \psi =\pi_- P\phi.
\end{equation}
Indeed if $\psi$ solves  \eqref{equa::mainpbm2}, then multiplying the equation by $\pi_-$ and assuming that the kernel of $\mathcal A$ is reduced to the nullspace (which is true except in the balanced case), one derives equation~\eqref{eq::newproblem}.
Such a problem admits a unique solution, given by the following theorem: 
\begin{thm}
\label{theo:invertibility}
The operator $\pi_- P \pi_-$ is invertible on $\mathcal R(\pi_-)$. Define $B_-$ the self-adjoint operator of $\mathcal{H}$ as  
\[B_-\phi=\pi_-(\pi_- P\pi_-)^{-1}\pi_-\phi.\]
Moreover it holds
\[ \Vert B_-\phi \Vert_\HH \le C\Vert \pi_-\phi \Vert_\HH \quad \text{ and } \quad
 \Vert B_-\phi \Vert_{1/2} \le C\Vert \pi_-\phi \Vert_{1/2}.\]
One can similarly define an operator $B_+$, obtained by replacing $\pi_-$ by $\pi_+$. 
\end{thm}

\FG{The result is proved} in \cite{fehrenbach2012generalized} for the full-Dirichlet case, that is $\Gamma_D=\partial \Omega$. \FG{The} proof can be adapted without major changes to the case in consideration. It is reproduced in Appendix \ref{sec:appendixpreuveth1} for the convenience of the reader.
Problem~\eqref{equa::mainpbm2} is then solved in the next theorem.

\begin{thm}\label{th:decomposition}
Let $\phi \in\mathcal H$, define $\Phi=(1,0)\in \mathcal{H}$ and consider problem \eqref{equa::mainpbm2}  of finding $\psi$  a solution of
\begin{equation}
\tag{\ref{equa::mainpbm2}}
\label{eq::pbm::total}
P\psi=P\phi \text{ and } \pi_+\psi=0.
\end{equation}
\begin{itemize}
\item If the constants are not controlled and $\int_\Omega h >  0$, there exists a solution  if and only if 
\[(\Phi-PB_-\Phi|\phi)_{\mathcal H}=0.\] 
In this case, the solution is unique and given \FG{by} $\psi=B_-P\phi$.
\item If the constants are not controlled and $\int_\Omega h =  0$ (balanced case), then $(\Phi-PB_-\Phi|\Phi)_{\mathcal H}\ne 0$ and there exists a unique solution given by
\[\psi=B_-P\phi+\frac{(\Phi-PB_-\Phi|\phi)_{\mathcal H}}{(\Phi-PB_-\Phi|\Phi)_{\mathcal H}}(B_-\Phi-\phi_0).\] 
\item In every other case then $\psi=B_-P\phi$ is the unique solution.
\end{itemize}
\end{thm}
The proof of this result is given in Appendix \ref{sec:proof:th:decomposition}.

\section{Resolution of the semi-infinite problem}
\label{sec:semi-infinite}
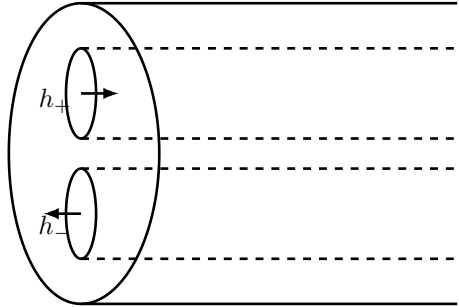
\begin{figure}
\begin{center}
\begin{tikzpicture}
\draw  [domain=0:360,scale=2,variable=\t,samples=200,line width=1pt] plot({0.5*cos(\t)+0.02*cos(\t)*cos(\t)}, {sin(\t)});
\draw [line width=1pt] (0,2)--(5,2);
\draw [line width=1pt] (0,-2)--(5,-2);

\draw  [domain=0:360,scale=2,variable=\t,samples=200,line width=1pt] plot({0.1*cos(\t)}, {0.4+0.3*sin(\t)});
\draw [line width=1pt,dashed] (0,1.4)--(5,1.4);
\draw [line width=1pt,dashed] (0,0.2)--(5,0.2);
\draw [line width=1pt,>=latex,->] (0,0.8)--(0.5,0.8);
\draw (0,0.7) node [left] {$h_+$}; 

\draw  [domain=0:360,scale=2,variable=\t,samples=200,line width=1pt] plot({0.1*cos(\t)}, {-0.4+0.3*sin(\t)});
\draw [line width=1pt,dashed] (0,-1.4)--(5,-1.4);
\draw [line width=1pt,dashed] (0,-0.2)--(5,-0.2);
\draw [line width=1pt,>=latex,->] (0,-0.8)--(-0.5,-0.8);
\draw (0,-0.7) node [below left] {$h_-$}; 
\end{tikzpicture}
\caption{Example of a semi-infinite cylinder, with two fluid domains.}\label{fig:semi-infini}
\end{center}
\end{figure}

In the semi-infinite problem, the equation is set on the cylinder $\Omega \times \RR^+$, see Figure \ref{fig:semi-infini}. The equation \eqref{eq::inside} becomes:
\[\partial_{zz} T +\Delta T +h\partial_z T= 0 \text{ on } \Omega \times \RR^+.\]
In this section we address different cases depending on the type of the Inlet/Oulet condition, namely either Dirichlet or Neumann.

In order to ensure uniqueness of the solution, we add the extra hypothesis that the temperature does not grow exponentially. We will say that the temperature has subexponential growth \FG{if and only if} for every $\lambda>0$, then $T(z)=o(e^{\lambda z})$ as $z$ goes to $+\infty$.

\subsection{Semi-infinite problem, Dirichlet Inlet/Outlet Condition}

We consider the Dirichlet \eqref{eq:ioboundary} condition:  
\begin{equation}\label{eq:conddirichlet}
T_{|z=0}=T_0 \quad \text{on } \Omega.
\end{equation}
 Denote  $\phi_{\mathcal D}=(T_0,0) \in \mathcal H$.

\begin{prop}
\label{prop:Semi-infinite-dirichlet}
Consider the Graetz problem \eqref{eq::inside} on the semi-infinite cylinder $\Omega\times[0,+\infty)$, with subexponential growth together with Dirichlet Inlet/Outlet condition \eqref{eq:conddirichlet}.
\begin{enumerate}
\item[a)] If the constants are controlled, then there exists a unique solution given by 
\[T=\sum_{i<0}e^{\lambda_iz} (B_-\phi_{\mathcal D}|\phi_i)_\HH U_i.\]
In this case the temperature at infinity is $0$.
\item[b)] If the constants are not controlled and $\int_\Omega h\ne 0$, then there exists a unique solution given by
\[T=\sum_{i<0}e^{\lambda_iz} (B_-(\phi_{\mathcal D}-T_\infty\Phi)|\phi_i)_\HH U_i+T_\infty,\]
where  $T_\infty$ is an arbitrary constant in the case $\int_\Omega h<0$ and $T_\infty=(\Phi-PB_-\Phi|\phi_{\mathcal D})_\HH(\Phi-PB_-\Phi|\Phi)_\HH^{-1}$ in the case $\int_\Omega h>0$.
In this case the temperature at infinity is the constant $T_\infty$. Note that if $\int_\Omega h>0$, the temperature at infinity is determined by $\phi_{\mathcal D}$ whereas in the case $\int_\Omega h<0$, it is a free parameter of the problem.
\item[c)] In the balanced case, the set of solutions is given by: 
\[T(z)=\sum_{i<0}(B_-(\phi_{\mathcal D}-c_1\Phi-c_2(\Delta_\sigma^{-1}h,0))|\phi_i)_\HH U_i e^{\lambda_i z} +c_1 +c_2(z+\Delta_\sigma^{-1}h),\]
where $c_2$ is an arbitrary constant and
\[c_1=(\Phi-PB_-\Phi|\phi_{\mathcal D}-c_2(\Delta_\sigma^{-1}h,0))_\HH(\Phi-PB_-\Phi|\Phi)^{-1}_\HH.\]
In this case the temperature at infinity has the linear growth rate $T(z)\simeq c_2 z+(c_2\Delta_\sigma^{-1} h+c_1)+o(z)$. If the temperature is not allowed to have a linear growth rate, then $c_2=0$ and the temperature at infinity is $c_1$ which is determined by the initial conditions.
\end{enumerate}
\end{prop}

\proof
We use the result of Section~\ref{sec:solevol} on the solution of the evolution equation, that is $T$ solves \eqref{eq::inside} \FG{if and only if} there exists $\psi \in \mathcal R(\mathcal A)$ and constants $c_1$ and $c_2$ such that:
\[T(z)=\sum_{i\in \ZZ^{*}}(\psi|\phi_i)_\HH U_i e^{\lambda_i z} +c_1 +c_2(z+\Delta_\sigma^{-1}h),\]
where $c_1=c_2=0$ if the constants are controlled and $c_2=0$ in the non-balanced case. The subexponential growth condition ensures that $\pi_+\psi=0$. The condition $(T(z=0),0)=\phi_{\mathcal D}$
 yields
 \begin{equation}
 \label{eq::to-be-solved:infty}
P\phi_{\mathcal D}=P\psi +c_1 \Phi +c_2(\Delta_\sigma^{-1}h,0).
\end{equation}
Using Theorem~\ref{th:decomposition} leads to distinguishing the following cases:\\
 a) If the constants are controlled then $c_1=c_2=0$ and the equation $P\phi_{\mathcal D}=P\psi$ with $\pi_+ \psi=0$ has the unique solution $\psi=B_-\phi_{\mathcal D}$. \\
b) In the non-balanced case, $c_2=0$ and $P(\phi_{\mathcal D}-c_1\Phi)=P \psi$ together with $\pi_+ \psi=0$ implies
 $\psi=B_-(\phi_{\mathcal D}-c_1\Phi)$ without any additionnal assumption in the case $\int_\Omega h<0$. In the case where $\int_\Omega h>0$ the compatibility condition is 
 \[(\Phi-PB_-\Phi|\phi_{\mathcal D}-c_1\Phi)=0, {\rm\ which\ gives\ } c_1=(\Phi-PB_-\Phi|\phi_{\mathcal D})(\Phi-PB_-\Phi|\Phi)^{-1}.\]
c) Finally, in the balanced case, let us fix an arbitrary value $c_2$. The conditions $\psi\in\mathcal R(\mathcal A)$ and $\pi_+\psi=0$ are equivalent to $\psi\in\mathcal R(\pi_-)$.
Denoting $\widetilde \Phi=(\Delta_\sigma^{-1}h,0)$ we re-cast \eqref{eq::to-be-solved:infty} into
\[P\phi_{\mathcal D}=P\psi +c_1\Phi+c_2\widetilde \Phi,\qquad\psi\in\mathcal R(\pi_-)\]
or equivalently
\[P(\phi_{\mathcal D}-c_2\widetilde \Phi)=P(\psi +c_1\phi_0),\qquad\psi\in\mathcal R(\pi_-).\]
In view of Theorem~\ref{th:decomposition} with $\phi=\phi_{\mathcal D}-c_2\widetilde \Phi$,  there is a unique solution to the above equation given by
\[ c_1=(\Phi-PB_-\Phi|\phi_{\mathcal D}-c_2\widetilde \Phi)_\HH(\Phi-PB_-\Phi|\Phi)^{-1}_\HH,\quad\psi=B_-(\phi_{\mathcal D}-c_1\Phi-c_2\widetilde \Phi).\] $\Box$

\subsection{Semi-infinite problem, Neumann Inlet/Outlet Condition}
We consider the Neumann \eqref{eq:ioboundary} condition:  
\begin{equation}\label{eq:condneu}
\partial_zT_{|z=0}=S_0.
\end{equation}
 Denote  $\phi_{\mathcal N}=(S_0,0) \in \mathcal H$.

\begin{prop}
Consider the Graetz problem \eqref{eq::inside} on the semi-infinite cylinder $\Omega\times[0,+\infty)$, with subexponential growth together with Neumann Inlet/Outlet condition \eqref{eq:condneu}.
\begin{enumerate}
\item[a)] If the constants are controlled, then there exists a unique solution given by 
\[T=\sum_{i<0}e^{\lambda_iz} (\mathcal A^{-1}B_-\phi_{\mathcal N}|\phi_i)_\HH U_i.\]
In this case the temperature at infinity is $0$.
\item[b)] If the constants are not controlled and $\int_\Omega h\ne 0$. If $\int_\Omega h>0$ there always exists a solution, if $\int_\Omega h<0$, there exists a solution \FG{if and only if} 
\[(\Phi-PB_-\Phi|\phi_{\mathcal N})=0.\]
When the solution exists, it is of the form
\[T=\sum_{i<0}e^{\lambda_iz} (\mathcal A^{-1}B_-(\phi_{\mathcal N}-T_\infty\Phi)|\phi_i)_\HH U_i+T_\infty,\]
where the temperature at infinity  $T_\infty$ is a free parameter of the problem.
\item[c)] In the balanced case, the set of solutions is given by 
\[T(z)=\sum_{i<0}(\mathcal A^{-1}B_-(\phi_{\mathcal N}+c_2\Phi)|\phi_i)_\HH U_i e^{\lambda_i z} +c_1 +c_2(z+\Delta_\sigma^{-1}h),\]
where $c_1$ is an arbitrary constant and $c_2$ is given by:
\[c_2=-(\Phi-PB_-\Phi|\phi_{\mathcal N})(\Phi-PB_-\Phi|\Phi)^{-1}.\]
In this case the temperature at infinity has the linear growth rate $T(z)\simeq c_2 z+(c_2\Delta_\sigma^{-1} h+c_1)+o(z)$.
\end{enumerate}
\end{prop}

\proof
We proceed as in the previous section. It follows from the \FG{result} of Section~\ref{sec:solevol} on the solution of the evolution equation that  $T$ solves \eqref{eq::inside} \FG{if and only if} there exists $\psi \in \mathcal R(\mathcal A)$ and constants $c_1$ and $c_2$ such that:
\[T(z)=\sum_{i\in \ZZ^{*}}(\psi|\phi_i)_\HH U_i e^{\lambda_i z} +c_1 +c_2(z+\Delta_\sigma^{-1}h),\]
where $c_1=c_2=0$ if the constants are controlled and $c_2=0$ in the non-balanced case. Derivating w.r.t. $z$ one finds
\begin{equation}
 \label{eq::to-be-solved:infty:Neumann}
P\phi_{\mathcal N}=P\mathcal A\psi +c_2 \Phi, \quad \pi_+\psi=0.
\end{equation}
Using Theorem~\ref{th:decomposition} leads to distinguishing the following cases:\\
a) If the constants are controlled then $c_1=c_2=0$ and equation \eqref{eq::to-be-solved:infty:Neumann} admits a unique solution $\mathcal A\psi=B_-\phi_{\mathcal N}$. The invertibility of $\mathcal A$ gives the result.\\
b) If the constants are not controlled then  $c_2=0$. If $\int_\Omega h > 0$ there is always a solution ${\mathcal A}\psi$ to equation \eqref{eq::to-be-solved:infty:Neumann} and the operator $\mathcal A$ is invertible, hence there exists a unique solution $\psi$ to \eqref{eq::to-be-solved:infty:Neumann} given by $\psi=\mathcal A^{-1}B_-\phi_{\mathcal N}$. The constant $c_1$ is then a free parameter of the problem.
If  $\int_\Omega h < 0$ then  the condition for equation \eqref{eq::to-be-solved:infty:Neumann} to admit a solution is:
 \[(\Phi-PB_-\Phi|\phi_{\mathcal N})=0.\]
 If this condition is met, by the invertibility of $\mathcal A$, $\psi=\mathcal A^{-1}B_-\phi_{\mathcal N}$ is the unique solution to \eqref{eq::to-be-solved:infty:Neumann} and $c_1$ is a free parameter of the problem.\\
 c) In the balanced case let $c_2$ be an arbitrary constant. It follows from Theorem~\ref{th:decomposition}  that $\mathcal A\psi$ satisfies \eqref{eq::to-be-solved:infty:Neumann} \FG{if and only if}
\begin{equation}
 \label{eq::to-be-solved:infty:Neumann:zero-debit}
\mathcal A\psi=B_-P\phi+\frac{(\Phi-PB_-\Phi|\phi)_{\mathcal H}}{(\Phi-PB_-\Phi|\Phi)_{\mathcal H}}(B_-\Phi-\phi_0)_\HH \quad{\rm where\ } \phi=\phi_{\mathcal N}+c_2\Phi.
\end{equation}
For $\psi$ to exist, the right hand side must belong to the range of $\mathcal A$, i.e. be orthogonal to $\phi_0$. Performing the scalar product of the left hand side of \eqref{eq::to-be-solved:infty:Neumann:zero-debit} with $\phi_0$ and recalling that the range of $B_-$ is orthogonal to $\mathcal K(\mathcal A)$ we obtain the following necessary condition:
\[(\Phi-PB_-\Phi|\phi_{\mathcal N}+c_2\Phi)_{\mathcal H}=0,\]
which is equivalent to: 
\[c_2=-(\Phi-PB_-\Phi|\phi_{\mathcal N})(\Phi-PB_-\Phi|\Phi)^{-1}.\]
Conversely, if the above condition is met, then equation \eqref{eq::to-be-solved:infty:Neumann:zero-debit} admits a unique inverse in $\mathcal R(\mathcal A)$ and $c_1$ is a free parameter of the problem. $\Box$

\section{Resolution of the problem in a finite domain}
\label{sec:finite}
We aim to solve Graetz equation in a domain of finite length $\Omega\times[-L,L]$:
\begin{equation}\label{eq:finitedirichlet}
\left\{\begin{array}{lr}
c\partial_{zz} T +\div \sigma \nabla T -h\partial_z T= 0 &\Omega\times[-L,L],\\
{\rm (LBC)} & (\partial\Omega)\times[-L,L],\\
{\rm Inlet/Outlet\ condition} & \Omega\times\{-L,L\},
\end{array}\right.
\end{equation}
where the I/O conditions can be of Neumann or Dirichlet type.
According to Section~\ref{sec:solevol}, the solutions may be sought in the form
$$T(z)=\sum_{i<0}(\psi|\phi_i)e^{\lambda_i(z+L)}U_i+\sum_{i>0}(\psi|\phi_i)e^{\lambda_i(z-L)}U_i +c_1 +c_2(z+\Delta_{\sigma}^{-1}h),$$
with $c_1=c_2=0$ if the constants are controlled and $c_2=0$ in the non-balanced case.

The unknowns in this equation are $(\psi|\phi_i)$ for $i<0$ and $i>0$, plus possibly (depending on the case) $c_1$ and $c_2$. Note that $\sum_{i<0}(\psi|\phi_i)U_i=P\pi_-\psi$, and therefore if $P\pi_-\psi$ is known \FG{it} suffices to decompose this vector on the basis of $L^2(\Omega)$ given by $(U_i)_{i<0}$ to obtain the desired coefficients for $i<0$. Similarly the coefficients $(\psi|\phi_i)$ for $i>0$ are obtained by considering the coefficients of $P\pi_+\psi$ on the basis composed of the $(U_i)_{i>0}$. Therefore the unknowns  to be determined are $P\pi_-\psi, P\pi_+\psi$ plus possibly $c_1$ and $c_2$.

Let $X$ be the vector composed of all the unknowns. Then satisfying the I/O conditions amounts to solving a linear system for $X$. In the rest of this section we detail the linear system in each case, but beforehand we  focus on a linear operator involved in the system.

\subsection{Study of the linear operator $M$}
\label{sec:studyM}
We define  and study  a linear operator that will be involved in the solution of the problem in a cylinder of length $2L$.
\begin{prop}
\label{defin:M}
Let $M_\pm$ be the operators from $\mathcal R(P)$ to $\mathcal R(P)$ and $M$ be given by:
\[M_\pm=P e^{\mp 2L\mathcal {A}}B_\pm \quad \text{ and } M=\begin{pmatrix} 0 & M_+ \\M_-& 0\end{pmatrix}.\]
Then 

\begin{enumerate}
\item[a)] there exists a constant $C$ such that
\[\Vert M\Vert \le Ce^{-2\lambda L}, \text{ where } \lambda=\min(\lambda_1,-\lambda_{-1}).\]
As a consequence $\Vert M\Vert<1$ for sufficiently large $L$.
\item[b)] If the constants are controlled, then for $L$ positive sufficiently small, $\Vert M^2\Vert<1$.
\item[c)] It follows that $Id+M$ is invertible on $\mathcal R(P)\times\mathcal R(P)$ for large $L$ and for small positive $L$.
\end{enumerate}
\end{prop}
\proof

a) Since $M_+=Pe^{-2L\mathcal A}B_+$ we have
$$\Vert M_+\Vert\le  \Vert B_+\Vert e^{-2L\lambda_1}.$$
A similar upper bound for $M_-$ gives the result.

b) Define
\[J(L)=\sup_{\Vert (\phi_1,\phi_2) \Vert=1}\Vert M^2(\phi_1,\phi_2) \Vert <1\]
Since $J(0)=1$, it is sufficient to prove that $J'(0)<0$.
Note that \[M^2=\begin{pmatrix} M_+M_- &0\\0&M_-M_+\end{pmatrix}.\]
Let us fix $\phi\in{\cal R}(P)$, and define $j(L)=\|M_+(L)M_-(L)\phi\|^2$. Then $j(0)=\|\phi\|^2$ and it remains to prove that $j'(0)\le -C\Vert \phi\Vert^2$ with a positive constant $C$ independent of $\phi$.
The derivative of $j$ is: 
\begin{eqnarray*}
j'(L) & = & (M_+'(L)M_-(L)\phi+M_+(L)M_-'(L)\phi|M_+(L)M_-(L)\phi)\\
 & = &(-2P\AA e^{-2L\AA}B_+Pe^{2L\AA}B_-\phi+2Pe^{-2L\AA}B_+P\AA e^{2L\AA}B_-\phi|Pe^{-2L\AA}B_+Pe^{2L\AA}B_-\phi)
\end{eqnarray*}
hence
\begin{eqnarray*}
j'(0)&=&-2(P\AA B_+\phi|\phi)+2(P\AA B_-\phi|\phi).
\end{eqnarray*}
But since $P\phi=\phi$, $P B_+P=P$ and $P\AA=\AA+P\AA P-\AA P$ we have
\begin{eqnarray*}
(P\AA B_+\phi|\phi)&=&(P\AA B_+\phi|B_+\phi)\\
&=&((\AA+P\AA P-\AA P)B_+\phi|B_+\phi)\\
&=&(\AA B_+\phi|B_+\phi)+(P\AA PB_+\phi|B_+\phi)-(\AA PB_+\phi|B_+\phi)\\
&=&(\AA B_+\phi|B_+\phi)+(\AA \phi|\phi)-(\AA \phi|B_+\phi).
\end{eqnarray*}
This proves that 
$$(P\AA B_+\phi|\phi)=(\AA \phi|B_+\phi)=\dfrac{1}{2}\left((\AA B_+\phi|B_+\phi)+(\AA \phi|\phi)\right).$$
Similarly we obtain that 
$$(P\AA B_-\phi|\phi)=\dfrac{1}{2}\left((\AA B_-\phi|B_-\phi)+(\AA \phi|\phi)\right).$$
As a summary we find that
$$
j'(0)=-(\AA B_+\phi|B_+\phi)+(\AA B_-\phi|B_-\phi)<\lambda_{-1}\|B_-\phi\|^2-\lambda_1\|B_+\phi\|^2<(\lambda_{-1}-\lambda_{1})\|\phi\|^2.
$$
c) The operator $Id+M$ is invertible for large $L$ by a).
Note that $M_\pm$ as endomorphism of $\mathcal R(P)$ are compact for $L>0$ and equal to identity for $L=0$. As a result $Id+M$ is invertible for small $L>0$ if and only if there is no eigenvector associated to the value $-1$. A sufficient condition for invertibility is then that $M^2$ does not admit $1$ as eigenvalue, which is proved in b) for $L$ sufficiently small. 
$\Box$

\subsection{The Dirichlet case}
The different cases for Dirichlet I/O condition are summarized in the following
\begin{prop}
The Dirichlet I/O condition $T|_{z=-L}=T_{-L}$ and $T|_{z=L}=T_{+L}$ are equivalent to the following linear system
$$ZX=b,$$
where $Z$, $X$ and $b$ are defined depending of the (LBC) and given in the table below.
\begin{tabular}{|c|c|c|c|}
\hline
& constants controlled & \multicolumn{2}{c|}{constants not controlled\hspace{2.8cm}\mbox{}}\\
&&unbalanced ($\int_\Omega h>0$)& balanced\\
\hline
$Z$ & $\begin{pmatrix} Id & M_+ \\M_-&Id \end{pmatrix}$ & $\begin{pmatrix}
Id & M_+ & \Phi\\
M_- & Id & \Phi\\
0 & M_+^\star u_-^T & (\Phi|u_-)
\end{pmatrix}$&
$\begin{pmatrix}
Id & M_+ & \Phi & -L\Phi+\tilde{\Phi}\\
M_- & Id & \Phi & L\Phi+\tilde{\Phi}\\
0 & M_+^\star u_-^T & (\Phi|u_-) & (-L\Phi+\tilde{\Phi}|u_-)\\
 M_-^\star u_+^T & 0 & (\Phi|u_+) & (L\Phi+\tilde{\Phi}|u_+)
\end{pmatrix}$\\
\hline
$X$&$\begin{pmatrix} P\pi_-\psi \\P\pi_+\psi \end{pmatrix}$&$\begin{pmatrix} P\pi_-\psi \\P\pi_+\psi\\c_1 \end{pmatrix}$&$\begin{pmatrix} P\pi_-\psi \\P\pi_+\psi\\c_1\\c_2 \end{pmatrix}$\\
\hline
$b$&$\begin{pmatrix}\phi_{-L} \\\phi_{+L} \end{pmatrix}$&$\begin{pmatrix}\phi_{-L} \\\phi_{+L} \\ (\phi_{-L}|u_-)\end{pmatrix}$&$\begin{pmatrix}\phi_{-L} \\\phi_{+L} \\ (\phi_{-L}|u_-) \\ (\phi_{L}|u_+) \end{pmatrix}$\\
\hline
\end{tabular}
where we recall that $\widetilde \Phi=(\Delta_\sigma^{-1}h,0)$ and we define $\phi_{\pm L}=(T_\pm,0)$,  $u_-=\Phi-PB_-\Phi$ and $u_+=\Phi-PB_+\Phi$.

Moreover, for 
sufficiently large $L$ this system is invertible. 
\end{prop}

Note 1: Thanks to the Lax-Milgram theorem in $3D$ (see Section~\ref{sec:lax-milgram}), we know before hand that there exists a unique solution to the system $ZX=b$.

\FG{Note 2}: In the case when the constants are not controlled and $\int_\Omega h<0$ it suffices to change the sign of $z$, or equivalently to remplace the $ _-$ by $ _+$.

\proof
The \eqref{eq:ioboundary} are equivalent to the following
\begin{equation}\label{eq:iodirichlet}
\left\{\begin{array}{l}
P\pi_-\psi = P\theta_- \qquad {\rm with\ }\theta_-=\phi_{-L}-e^{- 2L\mathcal {A}}\pi_+\psi-c_1\Phi-c_2(-L\Phi+\tilde{\Phi}),\\
P\pi_+\psi=  P\theta_+ \qquad {\rm with\ }\theta_+=\phi_{+L}-e^{ 2L\mathcal {A}}\pi_-\psi -c_1 \Phi-c_2(L\Phi+\tilde{\Phi}).
\end{array}\right.
\end{equation}
 Combining
\[Pe^{-2L\mathcal {A}}\pi_+\psi=Pe^{-2L\mathcal {A}}B_+P\pi_+\psi=M_+ P\pi_+\psi ,\]
and the similar version when the roles of $+$ and $-$ are interchanged with
equations~\eqref{eq:iodirichlet} we obtain the first two rows of the matrix $Z$.

a) when the constants are controlled, $c_1=c_2=0$ and \eqref{eq:iodirichlet}  reads $ZX=b$.

b) when the constants are not controlled and $\int_\Omega h>0$ then  $c_2=0$. Theorem \ref{th:decomposition} requires an additional compatibility condition to solve the first equation. This condition reads 
$$(\Phi-PB_-\Phi|\theta_-)=0,$$
which is the additional equation in the system $ZX=b$.

c) In the balanced case, after the change of variable $\widetilde \psi=\pi_-\psi$ the first equation in \eqref{eq:iodirichlet} $P\pi_-\psi=P\theta_-$ is equivalent to  
$$\begin{cases}
P\widetilde{\psi}=P\theta_-\\
\pi_+\widetilde\psi=0
\end{cases}\qquad{\rm and}\quad
(\widetilde\psi|\phi_0)=0.$$
Theorem \ref{th:decomposition} gives an explicit expression for the solution of the system on the left, and the condition on the right becomes $(\theta_-|\Phi-PB_-\Phi)=0$, which is the third row of the system $ZX=b$. The last row is  obtained using the second equation  in \eqref{eq:iodirichlet}.

When $L$ becomes large, $M_+$ and $M_-$ are exponentially small and in each case the matrix $Z$ is asymptotic to an invertible matrix.  The sole non-obvious case is the balanced case, where one can observe that the $2\times2$ lower right block 
is asymptotically equivalent to $\begin{pmatrix}
 (\Phi|u_-) & (-L\Phi|u_-)\\
 (\Phi|u_+) & (L\Phi|u_+)
\end{pmatrix}$
which has a determinant $2L(\Phi|u_-)(\Phi|u_+)\neq 0$. When $L$ is large, $Z$ can be rewritten as $Z=A+B$ with $B$ small and $A$ easily inverted. One can use a Neumann series strategy to solve $Zx=b$. 
$\Box$

\subsection{The Neumann Inlet/Outlet case}
\label{sec:finite:Neumann}
The different cases for Neumann I/O condition are summarized in the following
\begin{prop}
The Neumann I/O condition $\partial_zT|_{z=-L}=S_{-L}$ and $\partial_zT|_{z=L}=S_{+L}$ are equivalent to the following linear system
$$ZX=b,$$
where $Z$, $X$ and $b$ are defined depending of the (LBC) and given in the table below.
\begin{tabular}{|c|c|c|c|}
\hline
& constants controlled & \multicolumn{2}{c|}{constants not controlled\hspace{2.8cm}\mbox{}}\\
&&unbalanced ($\int_\Omega h>0$)& balanced\\
\hline
$Z$ & $\begin{pmatrix} Id & M_+ \\M_-&Id\end{pmatrix}$ & $\begin{pmatrix}
Id & M_+\\
M_- & Id\\
0 & M_+^\star u_-^T 
\end{pmatrix}$&
$\begin{pmatrix}
Id & M_+ & \Phi\\
M_- & Id & \Phi \\
0 & M_+^\star u_-^T & (\Phi|u_-) \\
 M_-^\star u_+^T & 0 & (\Phi|u_+) 
\end{pmatrix}$\\
\hline
$X$&$\begin{pmatrix} P{\mathcal A}\pi_-\psi \\P{\mathcal A}\pi_+\psi \end{pmatrix}$&$\begin{pmatrix} P{\mathcal A}\pi_-\psi \\P{\mathcal A}\pi_+\psi \end{pmatrix}$&$\begin{pmatrix} P{\mathcal A}\pi_-\psi \\P{\mathcal A}\pi_+\psi\\c_2 \end{pmatrix}$\\
\hline
$b$&$\begin{pmatrix}\phi_{-L} \\\phi_{+L} \end{pmatrix}$&$\begin{pmatrix}\phi_{-L} \\\phi_{+L} \\  (\phi_{-L}|u_-)\end{pmatrix}$&$\begin{pmatrix}\phi_{-L} \\\phi_{+L} \\ (\phi_{-L}|u_-) \\ (\phi_{L}|u_+) \end{pmatrix}$\\
\hline
\end{tabular}

where  we define $\phi_{\pm L}=(S_\pm,0)$,  $u_-=\Phi-PB_-\Phi$ and $u_+=\Phi-PB_+\Phi$.
\end{prop}

Note 1: when the constants are not controlled the value of $c_1$ is arbitrary. In these cases the linear systems are rectangular and the existence of the solution depends on a compatibility conditions that expresses that $b$ is in the range of $Z$.

Note 2: once the quantities $P{\mathcal A}\pi_\pm\psi$ are known, then the $({\mathcal A}\psi|\phi_i)$ for $i>0$ and $i<0$ can be computed as explained above, and $(\psi|\phi_i)$ is obtained by dividing by $\lambda_i$.

\proof
The \eqref{eq:ioboundary}  are equivalent to the following
\begin{equation}\label{eq:ioneumann}
\left\{\begin{array}{l}
P{\mathcal A}\pi_-\psi = P\theta_- \qquad {\rm with\ }\theta_-=\phi_{-L}-e^{- 2L\mathcal {A}}{\mathcal A}\pi_+\psi-c_1\Phi-c_2(-L\Phi+\tilde{\Phi}),\\
P{\mathcal A}\pi_+\psi=  P\theta_+ \qquad {\rm with\ }\theta_+=\phi_{+L}-e^{ 2L\mathcal {A}}{\mathcal A}\pi_-\psi -c_1 \Phi-c_2(L\Phi+\tilde{\Phi}).
\end{array}\right.
\end{equation}
A discussion similar to the Dirichlet case leads to the result.$\Box$

\section{Numerical tests}
\label{sec:numerical}
\subsection{First test case: a domain of finite length}
\begin{figure}[!hbt]
\center
\includegraphics[width=0.5\textwidth]{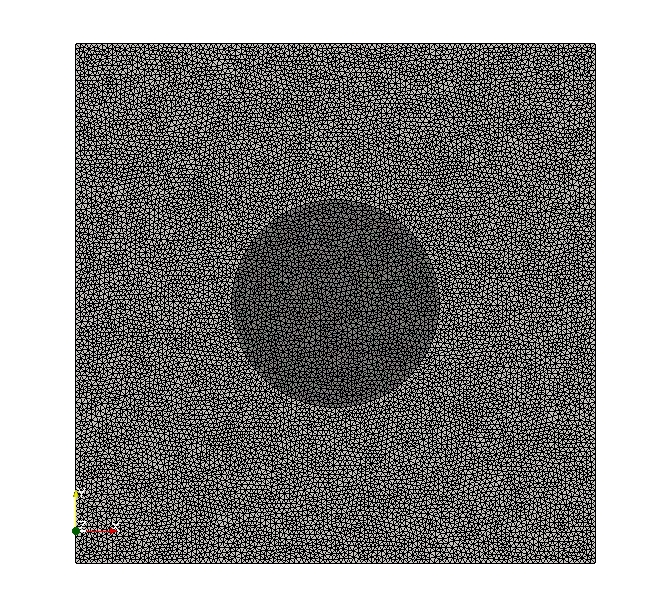}
\caption{The mesh for the first test case is composed of 13589 vertices and 26776 triangles. The solid domain is in white and the fluid domain in grey.}\label{fig:maillage}
\end{figure}

The section of the domain of the first test-case is the square $\Omega=[-5,5]^2$ with a circular fluid subdomain of radius $2$ centered at the origin.
The velocity and eigenvalues of the operator $\mathcal A$ are computed with $P1$ finite element methods on the mesh of Figure~\ref{fig:maillage} .
The velocity has a parabolic profile (Poiseuille flow) with prescribed total flow $Q \in\{1,10,100,1000\}$. The lateral boundary conditions are of Robin type with parameter $a$. The thermal conductivities are equal to $c=\sigma=1$. In total $100$ eigenvalues / eigenvectors of $\mathcal A$ are computed.

\pgfplotstableread[col sep=comma]{output4.dat}{\DTN}
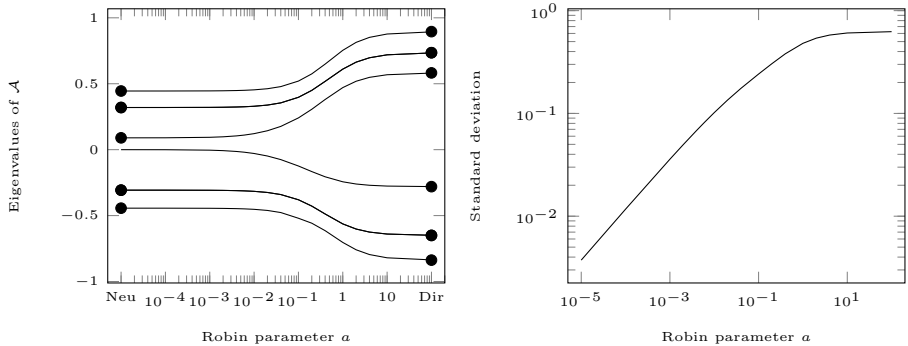
\begin{figure}[!hbt]
\begin{center}
\begin{tikzpicture}[every node/.append style={font=\tiny}]
\begin{axis}[width=0.5\textwidth,xmode=log, xmin=0.000005, xmax= 200, xlabel={Robin parameter $a$}, ylabel={Eigenvalues of $\mathcal A$},
xtick={1e-5,1e-4,1e-3,1e-2,1e-1,1,10,100},
xticklabels={Neu,$10^{-4}$,$10^{-3}$,$10^{-2}$,$10^{-1}$,$1$,$10$,Dir},
]
\addplot [color=black] table[y={1}]{\DTN};
\addplot [color=black] table[y={2}]{\DTN};
\addplot [color=black] table[y={3}]{\DTN};
\addplot [color=black] table[y={4}]{\DTN};
\addplot [color=black] table[y={6}]{\DTN};
\addplot [color=black] table[y={7}]{\DTN};
\addplot [color=black] table[y={8}]{\DTN};
\addplot [color=black] table[y={9}]{\DTN};
\addplot [color=black,only marks,mark=*]
table{
x    y
0.00001 0.089669419412
0.00001 0.319734873043
0.00001 0.319734883535
0.00001 0.444872847387
0.00001 -0.307047772647
0.00001 -0.30704782975
0.00001 -0.443666697998
};
\addplot+[color=black,only marks,mark=*,mark options={solid},]
table{
x    y
100 0.582352943623
100 0.735370979774
100 0.735371039889
100 0.895273148282
100 -0.280497253995
100 -0.650375910117
100 -0.650374654948
100 -0.837718141943
};
\end{axis}
\end{tikzpicture}
\begin{tikzpicture}[every node/.append style={font=\tiny}]
\begin{axis}[width=0.5\textwidth,xmode=log,ymode=log, xmin=0.000005, xmax= 200, xlabel={Robin parameter $a$}, ylabel={Standard deviation}]
\addplot[color=black,]
table{
x    y
0.00001 0.00373749502274
0.0001 0.0117038116388
0.001 0.0358416093204
0.002 0.0496789284584
0.004 0.068249033315
0.006 0.0817275411297
0.01 0.101802148746
0.02 0.134971171406
0.04 0.175337643049
0.1 0.241594062815
0.2 0.304918732337
0.4 0.380107473538
1. 0.48080436013
2 0.539965993888
4 0.578811806376
10. 0.606202319394
100. 0.624413952814
};
\end{axis}
\end{tikzpicture}
\end{center}
\caption{On the left: Evolution of the eigenvalues of $\mathcal A$ of smallest magnitude for varying parameter of the Robin lateral boundary condition. Eigenvalues for the Neumann (resp.  Dirichlet) boundary conditions are shown as bullet on the left (resp. right) of the curves. On the right: relative $L^2$ difference between the eigenvector with largest negative eigenvalue and its mean as the Robin parameter varies.}
\label{fig:Robinvaries}
\end{figure}
We first set $Q=10$ and vary the Robin parameter $a$. When $a=0$, one retrieves the Neumann case 
and when $a=+\infty$, one retrieves the Dirichlet case. 
In order to emphasize this fact we plot in Figure~\ref{fig:Robinvaries} (left) the eigenvalues of smallest magnitude for different values of $a$. We also plot with dots the eigenvalues associated to the Neumann problem (on the left of the curves) and the one associated to the Dirichlet case (on the right of the curves).
The smooth transition \FG{from} Neumann to Dirichlet as the Robin parameter varies is striking except from the fact that there exists an eigenvalue that goes to zero as $a$ goes to zero even if the Neumann problem does not have zero \FG{as} eigenvalue. 
We claim that this behavior is consistent with theory. First $0$ is not an eigenvalue of the Neumann case since the total flow is non-zero (hence we are not in a balanced case even if the constants are not controlled). Second, we remark that the zero eigenvalue is the limit of a negative eigenvalue. Remember from Proposition~\ref{prop:Semi-infinite-dirichlet} that it is always possible to decompose a temperature field on the set of negative eigenvectors in the Robin case (part a) of the proposition), but for the Neumann case it is necessary to add \FG{a} constant (part b) of the proposition). In other words, in the Neumann case the constant must be added to the negative eigenvectors to obtain a Hilbert basis of $H$, while the set of positive eigenvectors form a Hilbert basis on their own. This explains why the constant emerges as the limit of a negative eigenvector, see Figure~\ref{fig:Robinvaries} (right) where the convergence of the eigenvector to the constant is numerically demonstrated.

\medskip

\pgfplotstableread[col sep=comma]{output.dat}{\EvolVp}
\pgfplotstableread[col sep=comma]{output2.dat}{\EvolVpDeux}
\pgfplotstableread[col sep=comma]{output3.dat}{\EvolVpTrois}

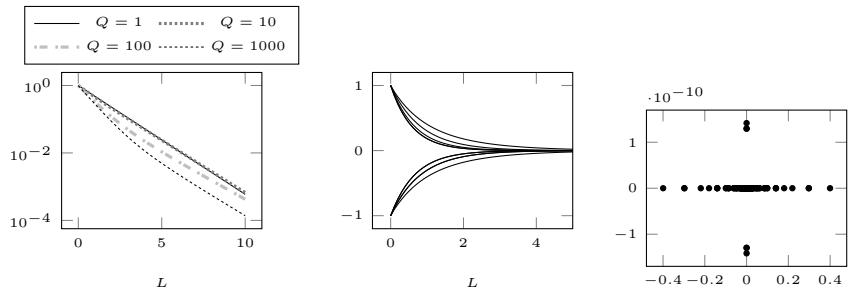
\begin{figure}[!hbt]
\begin{tabular}{ccc}
\begin{tikzpicture}[every node/.append style={font=\tiny}]
\begin{axis}[ymode=log,width=0.35\textwidth,legend style={ at={(0.5,1.05)}, anchor=south, legend columns=2},xlabel=$L$]
\addplot [color=black]
          table[y={1}]{\EvolVp};
\addplot [color=gray,dash pattern=on 1pt off 1pt,very thick]
          table[y={10}]{\EvolVp};
\addplot [color=lightgray,dash pattern=on 1pt off 1pt on 3pt off 3pt,very thick]
          table[y={100}]{\EvolVp};
\addplot [color=black,dash pattern=on 1pt off 1pt]
          table[y={1000}]{\EvolVp};
\legend{$Q=1$,$Q=10$,$Q=100$,$Q=1000$}
\end{axis}
\end{tikzpicture}
&
\begin{tikzpicture}[every node/.append style={font=\tiny}]
\begin{axis}[width=0.35\textwidth,xmax=5,xlabel=$L$]
\addplot [color=black] table[y={1}]{\EvolVpDeux};
\addplot [color=black] table[y={2}]{\EvolVpDeux};
\addplot [color=black] table[y={3}]{\EvolVpDeux};
\addplot [color=black] table[y={4}]{\EvolVpDeux};
\addplot [color=black] table[y={4}]{\EvolVpDeux};
\addplot [color=black] table[y={196}]{\EvolVpDeux};
\addplot [color=black] table[y={197}]{\EvolVpDeux};
\addplot [color=black] table[y={196}]{\EvolVpDeux};
\addplot [color=black] table[y={199}]{\EvolVpDeux};
\addplot [color=black] table[y={200}]{\EvolVpDeux};
\end{axis}
\end{tikzpicture}
&
\begin{tikzpicture}[every node/.append style={font=\tiny}]
\begin{axis}[width=0.35\textwidth]
\addplot [color=black,mark=*,only marks,mark size=1pt] table {\EvolVpTrois};
\end{axis}
\end{tikzpicture}
\end{tabular}
\caption{Eigenvalues of $M$ for a Robin test case with $a=1$. Left: evolution of the spectral radius for different values of the total flow $Q$. Center: the evolution of the five largest positive and first smallest negative eigenvalue for a total flow of $Q=20$. Right: eigenvalues in the complex plane for $L=1$ and $Q=20$.}
\label{Fig:Suivi-VP}
\end{figure}
In a second parametric study, we fix $a=1$ and we let both $Q$ and $L$ vary. First we plot the spectral radius of the matrix $M$ defined in Proposition~\ref{defin:M} versus the exchanger length $L$ for the different values of the total flow $Q$ in Figure~\ref{Fig:Suivi-VP} (left). Figure ~\ref{Fig:Suivi-VP} (center) shows the evolution of the 5 smallest positive and 5 largest negative eigenvalues of $M$ for a fixed total flow $Q=20$.
This test-case shows that, apart from the case $L=0$, the spectral radius of the matrix $M$ is always smaller than one, so that the matrix $Id+M$ is indeed always invertible. 
The exponential decrease for large $L$ and the decrease at the origin follows from Proposition~\ref{defin:M}. Moreover, since the spectral radius of $M$ is strictly smaller than one, a Neumann series strategy to solve 
$$(Id+M)^{-1}b=\sum_k (-M)^kb$$
is legitimate.
In Figure~\ref{Fig:Suivi-VP} (right), the whole spectrum of $M$ is shown in the complex plane. Although the spectrum seems real, we do not have mathematical proof of this fact.

\subsection{Second test-case: a periodic exchanger}

\begin{figure}[!hbt]
\begin{center}
\includegraphics[width=0.3\textwidth]{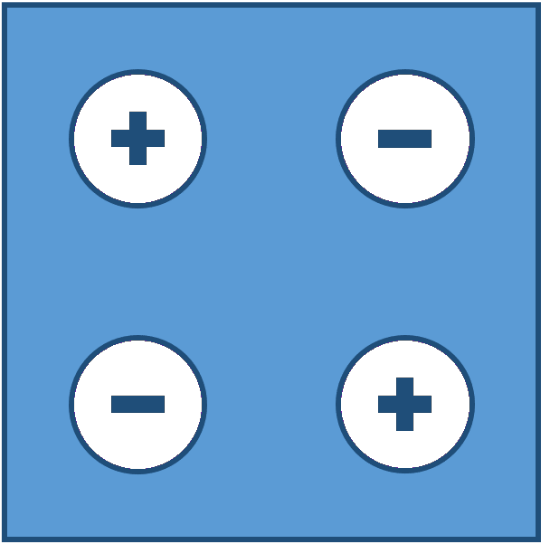}
\includegraphics[width=0.6\textwidth]{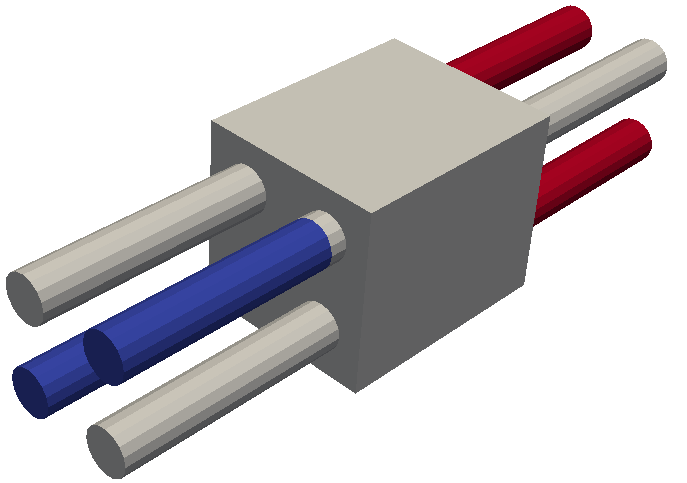}
\end{center}
\caption{Geometry of the periodic exchanger. On the left, a cut inside the exchanger with the sign of the fluid velocities. On the right, a $3D$ representation of the exchanger. The tubes where the temperature is set at $\infty$ are colored accordingly to their temperature.}\label{fig:explainexch}
\end{figure}

\begin{figure}[!hbt]
\begin{center}
\includegraphics[width=0.45\textwidth]{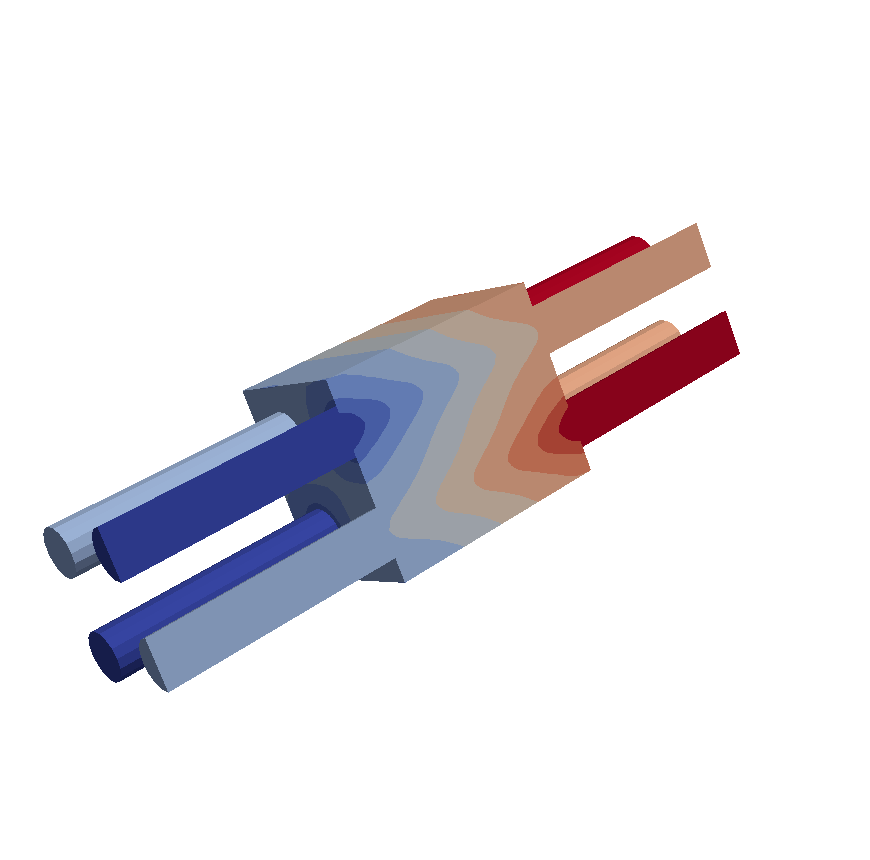}
\includegraphics[width=0.45\textwidth]{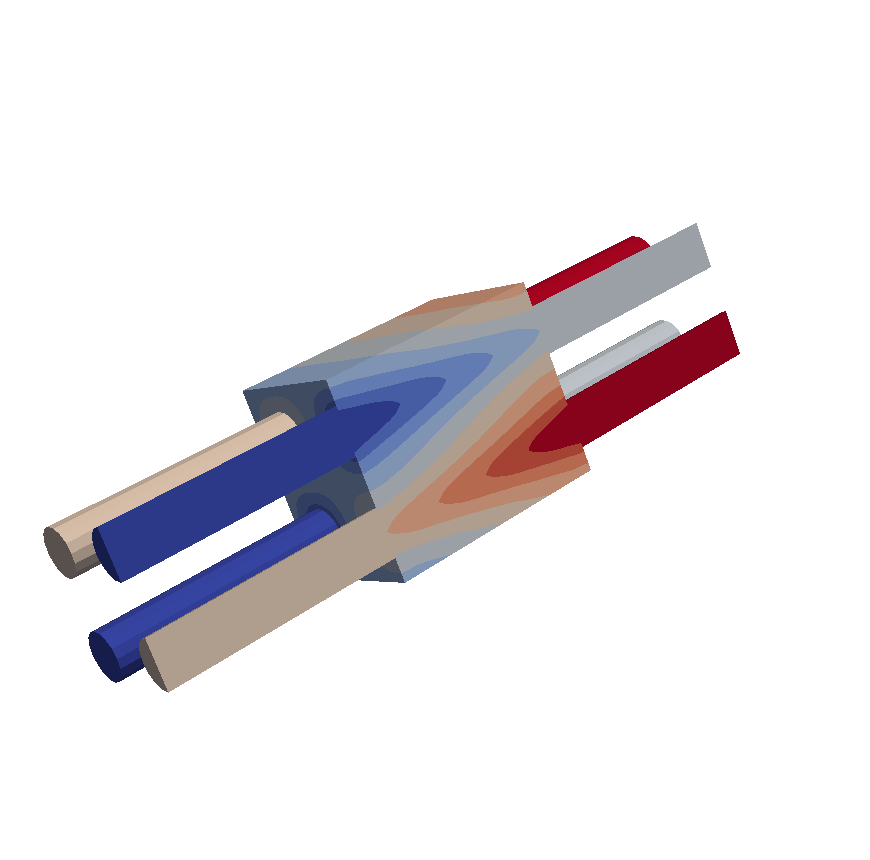}
\includegraphics[width=0.45\textwidth]{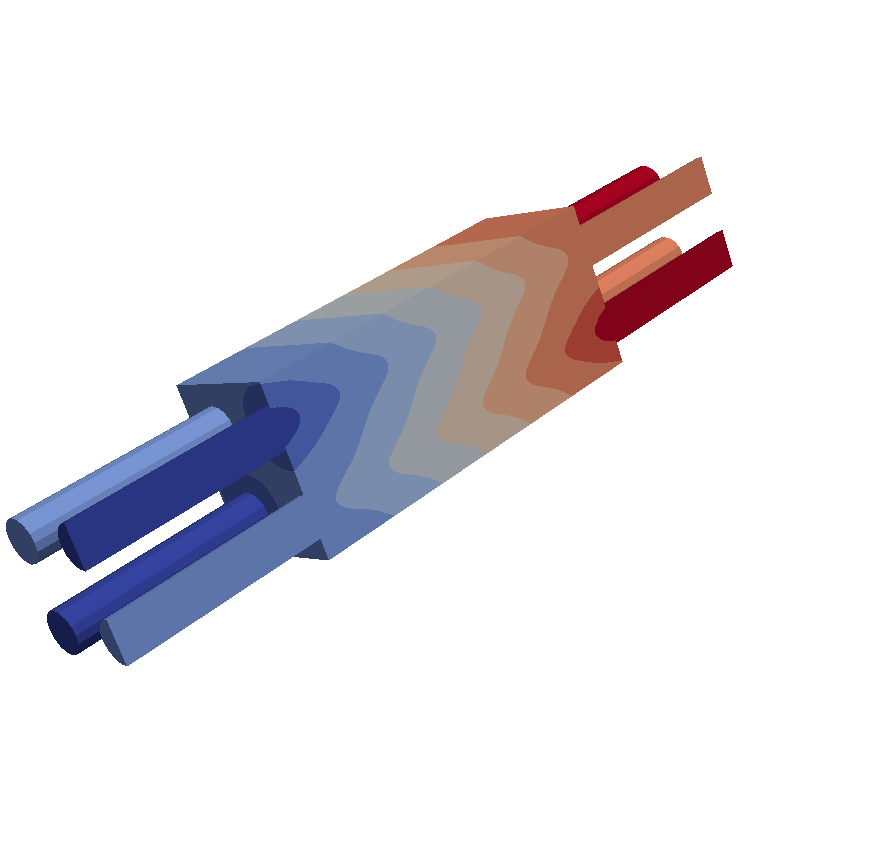}
\includegraphics[width=0.45\textwidth]{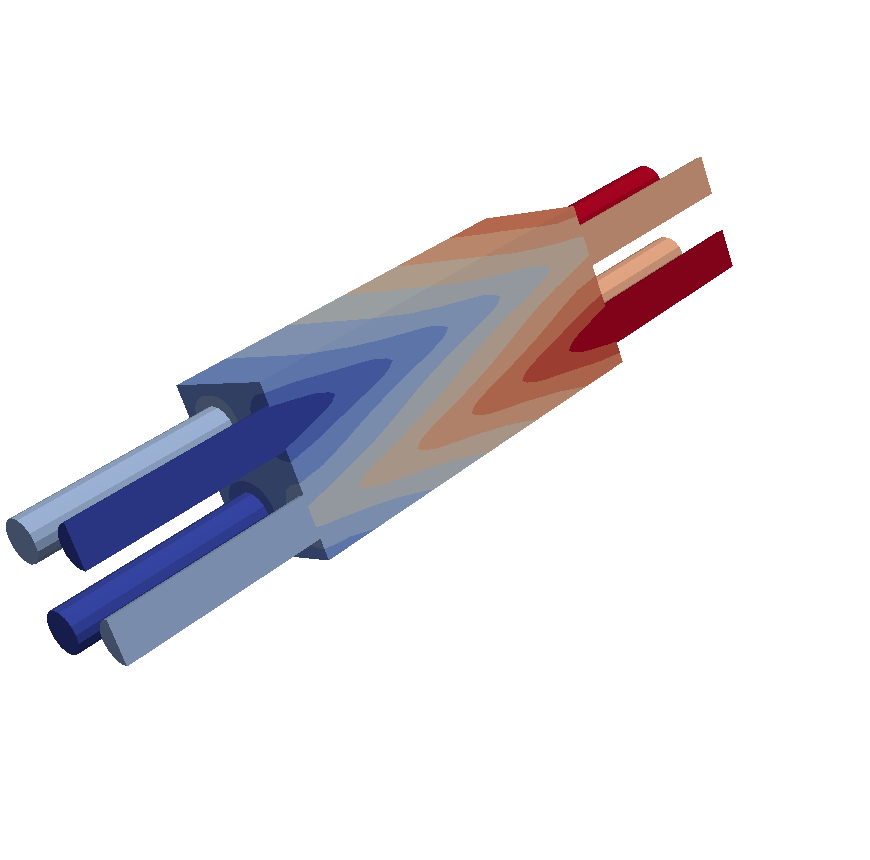}
\end{center}
\caption{Four different solutions of the periodic exchanger, with different length and total flow. The length of the exchanger is set to $L=10$ on top and $L=20$ on bottom. The total flow is set to   to $Q=10$ on left and $Q=30$ on right.}\label{fig:sol3d}
\end{figure}

The second test case consists of a heat exchanger with periodic boundary conditions. The whole device consists of one solid exchanger through which pass four tubes containing fluids. A cut along the middle of the exchanger is shown in Figure~\ref{fig:explainexch} (left) where the sign of the velocity of the fluid in the inner tubes is displayed. The fluids are assumed to obey a Poiseuille flow, the velocities are then quadratic in the radial coordinates of their corresponding tubes. The length of the exchanger is denoted $L$, the section of the exchanger is the square $[-4,4]^2$, the radii of the inner tubes \FG{are} fixed to $1$ and the distance of the center of the inner tubes to the center of the exchanger is $\sqrt{2^2+2^2}$.
The conductivity in both the fluid and solid part is set to $1$. The temperature is fixed for the four tubes with incoming flow (two at each side) on the exchanger, the warm temperature is set to $+1$ and the cold temperature to $-1$, see Figure~\ref{fig:explainexch} (left).
In what follows, $Q$ denotes the total flow of fluid in one tube. We glue together the different Graetz problems using the methodology developped in~\cite{fehrenbach2012generalized}.

\begin{figure}[!hbt]
\center
\begin{tikzpicture}
\begin{axis}[
width=0.45\textwidth,
view/h=-60,
view/v=60,
colorbar horizontal,
xlabel=$L$, ylabel=$Q$,
]
\addplot3[surf,mesh/ordering=y varies]
table {efficiency.dat};
\end{axis}
\end{tikzpicture}
\begin{tikzpicture}
\begin{axis}[
width=0.45\textwidth,
view/h=-60,
view/v=60,
colorbar horizontal,
xlabel=$L$, ylabel=$Q$,
]
\addplot3[surf,mesh/ordering=y varies]
table {exchange.dat};
\end{axis}
\end{tikzpicture}
\caption{Values of the efficiency (left) and the exchange (right) for different values of the length $L$ and the total flow $Q$ of the periodic exchanger. $L$ ranges from \FG{$0.5$} to $13$ and $Q$ ranges from $1$ to $30$. Each direction has been sampled $50$ times for a total of $2500$ exchanger computations.}
\label{fig::echange and efficiency}
\end{figure}
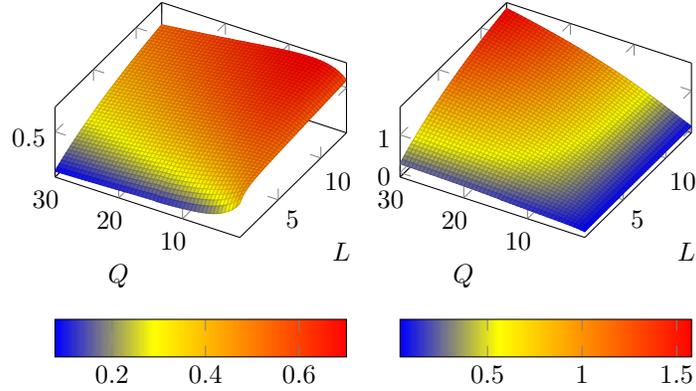

In Figure~\ref{fig:sol3d}, four solutions are shown for different values of the length $L$ and the flow $Q$.

Figure ~\ref{fig::echange and efficiency} displays the efficiency and the total exchange for different values of $Q$ and $L$. For a tube containing fluid whose velocity is positive (resp. negative), the temperature at $-\infty$ (resp. $+\infty$) is set to $1$ (resp $-1$), the efficiency of the exchanger is then defined by $-T_{+\infty}/T_{-\infty}$ (resp $-T_{-\infty}/T_{+\infty}$), where $T_{\pm \infty}$ is the temperature at infinity. This efficiency is between $-1$ and $1$. The exchange is simply the total amount of heat exchanged and is equal to $Q$ times the efficiency.
The aim of this test-case is to document the fact that our method is able to deal with any boundary conditions and type of exchanger. It is well suited for parametric studies.

\section{Conclusion}
In the present work we have proposed a general framework dedicated to the resolution of the generalized Graetz problem in arbitrary geometry, involving any type of boundary \FG{conditions}. The main novelty is the introduction of insulating boundary condition (Neumann or periodic) that allow to model realistic heat exchangers. Our study highlighted a special case that we call the balanced case, when $\int_\Omega h=0$ (together with Neumann or periodic boundary condition) where the solution is different than in the general case. We have also proposed a number of numerical illustration in various test cases.

\appendix
\section{Technical lemmas}
\label{sec:decomp:notation}
We prove here results that will be used in the sequel.


\begin{lem}
For each $\phi=(u,s) \in {\mathcal D(\mathcal A)},\tilde \phi=(\tilde u, \tilde s) \in {\mathcal H}$, we have
\begin{eqnarray}
(Id-P){\mathcal A}(Id-P)\phi&=&0 \label{sp.1a}. \\
 (P\mathcal A P\phi|\tilde \phi)_{\mathcal H}&=&\int_{\Omega} hu\tilde u \le \Vert h\Vert_{L^\infty(\Omega)} \Vert P\phi\Vert_{\mathcal H}\Vert P\tilde \phi\Vert_{\mathcal H} \label{sp.1b}.
\end{eqnarray}
\end{lem}
\proof
This results from elementary calculations using  the definition of $\mathcal A$ \eqref{eq:def:A}, and  the definition of $P$ \eqref{eq:def:P}.
$\Box$

\begin{lem}
\label{lemma::properties::Amoins1}
 Let $\Phi=(1,0)\in \mathcal H$. Let $\phi\in \mathcal D(A^{-1/2})$ such that $P\phi=\phi$, then
\begin{eqnarray*}
({\mathcal A}^{-1}\phi |\phi)_{\mathcal H}&=&\left(\int_\Omega h\right)^{-1}(\phi|\Phi)_{\mathcal H}^2  \text{ in the non-balanced case} \\
({\mathcal A}^{-1}\phi |\phi)_{\mathcal H}&=& 0 \text{ in the ``balanced'' or ``constant controlled'' case}
\end{eqnarray*}
\end{lem}
\proof
The expression of $\mathcal A^{-1}$ is given in section \ref{sec:operatorA} for the various cases.
Let $\phi \in \mathcal D( \mathcal A^{-1/2})$ such that $P\phi=\phi$, hence there exists $u\in L^2(\Omega)$ such that $\phi=(u,0)$.
 If the constants are controlled then ${\mathcal A}^{-1}\phi=(0,\Delta_\sigma^{-1}(-cu))$ and $({\mathcal A}^{-1}\phi|\phi)_\HH=0$. 
 If the constants are not controlled and $\int_\Omega h\ne 0$, then 
 \[{\mathcal A}^{-1}\phi=(k,\Delta_\sigma^{-1}(-cu+hk)) \text{ with } k\int_\Omega h=\int_\Omega cu,\]
 hence \[({\mathcal A}^{-1}\phi|\phi)_\HH=\int_\Omega cku=(\int_\Omega h)^{-1}(\int_\Omega cu)^2.\]
Finally, in the balanced case, since $\phi \in   \mathcal D( \mathcal A^{-1/2})$, then $(\phi|\phi_0)_\HH=0$ and ${\mathcal A}^{-1}\phi=
(0,\Delta_\sigma^{-1}(-cu)) + k\phi_0$ exists and $({\mathcal A}^{-1}\phi|\phi)_\HH=0$.$\Box$

\begin{lem}\label{lem:lem5}
In the balanced case, 
$$(\Phi-PB_-\Phi| \Phi)_\HH\ne 0.$$
\end{lem}
\proof
 Suppose the contrary and set $\theta=\Phi-PB_-\Phi$, we have $P\theta= \theta$ and  by definition of $B_-$, we have $\pi_- \theta=0$. Moreover, we have
\begin{eqnarray*}
( \theta|\phi_0)_\HH&=&(\Phi-PB_-\Phi|\phi_0)_\HH=(P\Phi-PB_-\Phi|\phi_0)_\HH=(\Phi-B_-\Phi|P\phi_0)_\HH\\
&=&(\Phi-B_-\Phi| P\Phi)_\HH=(\Phi-PB_-\Phi| \Phi)_\HH=0
\end{eqnarray*}
Lemma~\ref{lemma::properties::Amoins1} ensures that $\mathcal A^{-1}\theta$ exists and that
  \[(\mathcal A^{-1} \theta| \theta)_\HH=0\]
  Since $ \theta$ belongs to $\mathcal R(\pi_+)$ and  all the eigenvalues of $\mathcal A$ are  positive on this space, this implies that $ \theta=0$ and then  $\Phi=PB_-\Phi$. Hence there exists $s$ such that $B_-\Phi=\begin{pmatrix}1, s\end{pmatrix}$ and 
\[({\mathcal A}B_-\Phi|B_-\Phi)_\HH= (\begin{pmatrix}h-\Delta_\sigma s,0\end{pmatrix} |\begin{pmatrix}1, s\end{pmatrix})_\HH=0.\]
But $B_-\Phi$ belongs to $\mathcal R(\pi_-)$ and since all the eigenvalues of $\mathcal A$ are  negative on $\mathcal R(\pi_-)$, $({\mathcal A}B_-\Phi|B_-\Phi)_\HH=0$ implies that $B_-\Phi=0$ which is in violation of $\Phi=PB_-\Phi$. Hence $(\Phi-B_-\Phi| \Phi)_\HH\ne 0$.$\Box$

\section{Proof of Theorem~\ref{theo:invertibility}}
\label{sec:appendixpreuveth1}
Let $M\in \NN$  and denote for short $\pi=\pi_{\integer{-M}{-1}}$. The operator $\pi P\pi$ is a symmetric operator on a finite-dimensional space, hence it is diagonalisable in an orthonormal basis. The first step is to prove that this operator is definite positive with a lower bound on its \FG{eigenvalues} that is independent of $M$. Let $\rho$ be an eigenvalue of $\pi P \pi $ and $\vv$ an associated normalized eigenvector: 
$\pi P \pi \vv=\rho \vv$, $(\vv|\vv)_{\mathcal H}=1$ and $\pi \vv=\vv$. Since
\[\rho=(\pi P\pi \vv|\vv)_\HH=(P\pi \vv|\pi \vv)_\HH=(P\pi \vv|P\pi \vv)_\HH=\Vert P\vv\Vert_\HH^2 \le \Vert \vv\Vert_\HH^2=1,\]
then $0\le \rho\le 1$.
Using \eqref{sp.1b} gives
\[|(P{\mathcal A}P\vv|\vv)_\HH|\le \Vert h\Vert_{L^\infty(\Omega)}\Vert P\vv\Vert_\HH^2.\]
It follows from \eqref{sp.1a} that $((Id-P){\mathcal A}(Id-P)\vv|\vv)_\HH=0$ and $\pi {\mathcal A}={\mathcal A}\pi$, we have
\[(P{\mathcal A}P\vv|\vv)_\HH=(2\rho -1)({\mathcal A}\vv|\vv)_\HH\]
Since $|({\mathcal A}\vv|\vv)_\HH|=|\sum_{i\in I} \lambda_i \frac{(\vv|\phi_i)_\HH^2}{\Vert \phi_i\Vert_\HH^2}|\ge |\lambda_{-1}|\Vert \vv\Vert_\HH^{2}= |\lambda_{-1}|$, we have 
\begin{equation}
\label{eq::bnd.gamma}
|\lambda_{-1}(2\rho -1)|\le \Vert h\Vert_{L^\infty(\Omega)}\Vert P\vv\Vert_\HH^2=\Vert h\Vert_{L^\infty(\Omega)}\rho. 
\end{equation}
This \FG{in turn} implies that $\rho \ge \dfrac{|\lambda_{-1}|}{2|\lambda_{-1}|+\Vert h\Vert_{L^{\infty}(\Omega)}}$, hence there exists $C$ independent of $M$ such that
\begin{equation}
\label{eq::coercivity}
(\pi P \pi \phi |\phi)_\HH\ge C\Vert \pi \phi\Vert_\HH \quad \forall \phi \in \mathcal{H(A)}.
\end{equation}
Since $\pi_-\phi$ is the strong $\mathcal H$-limit of $\pi\phi$ as $M$ goes to infinity and the constant $C$ does not depend on $M$. Passing to the limit, we recover equation~\eqref{eq::coercivity} with $\pi$ replaced by $\pi_-$. The Lax-Milgram theorem applies and $\pi_- P \pi_-$ is a bijection \FG{from} $\mathcal R(\pi_-)$ onto $\mathcal R(\pi_-)$ with a continuous inverse bounded by a constant in $\mathcal H$-norm. 

We turn our interest to the bound in ${1/2}$ norm of $B_-$. Let $\phi\in \mathcal D(\mathcal A^{1/2})$, for any $M\in \NN$ denote $\pi=\pi_{\integer{-M}{-1}}$, and let  $\vv=\pi B_-\phi$. We have $\pi \vv=\vv$ and $\vv \in \mathcal D(\mathcal A)$. Recalling~\eqref{sp.1a} and $\pi \mathcal A={\mathcal A}\pi$, , we have 
\[
(P{\mathcal A}P\vv|\vv)_\HH=(({\mathcal A}P+P{\mathcal A}-{\mathcal A})\vv|\vv)_\HH=2(P\vv,{\mathcal A}\vv)_\HH-({\mathcal A}\vv,\vv)_\HH.
\]
Hence, since $\pi \vv=\vv$ and $\pi$ is a projection on negative eigenvalues of $\mathcal A$ only, then $\Vert\vv\Vert_{1/2}^{2}=-({\mathcal A}\vv|\vv)_\HH$ and 
\begin{equation}
\label{eq:bnd1/2}
\Vert \vv \Vert_{1/2}^2=(P{\mathcal A}P\vv|\vv)_\HH-2(P\vv|{\mathcal A}\vv)_\HH\le \Vert h\Vert_{L^\infty(\Omega)}\Vert \vv\Vert_\HH^2+2\Vert \pi P\vv \Vert_{1/2}\Vert \vv \Vert_{1/2}
\end{equation}
Using the bound on the $\mathcal H$-norm of $B_-$, we have 
\begin{equation}
\label{eq:bnd1/22}
\Vert \vv\Vert_\HH=\Vert \pi B_-\phi\Vert_\HH \le C\Vert \pi_- \phi\Vert_\HH\le C\Vert \pi_- \phi\Vert_{1/2}.
\end{equation}

We infer from \eqref{eq:bnd1/2} and \eqref{eq:bnd1/22} that $\Vert \vv \Vert_{1/2}\le C(\Vert \pi_- \phi \Vert_{1/2}+\Vert \pi P\vv \Vert_{1/2})$. We let $M$ go to infinity, then $\pi P\vv=\pi P\pi B_-\phi$ goes to $\pi_-\phi$ and $\vv$ goes to $B_-\phi$, we obtain:
\[\Vert B_- \phi \Vert_{1/2}\le C\Vert \pi_-\phi \Vert_{1/2},\]
which finishes the proof.

\section{Proof of Theorem~\ref{th:decomposition}}
\label{sec:proof:th:decomposition}

{\em First case: $ \mathcal K(\mathcal A)=\{0\}$, i.e. every case but the balanced case}

In this case the condition $\pi_+\psi=0$ is then equivalent to $\psi=\pi_-\psi$.

After multiplication of \eqref{eq::pbm::total} by $B_-\pi_-$, one obtains the following necessary condition for \eqref{eq::pbm::total} to hold, which proves uniqueness: 
$$\psi=B_-P\phi.$$ 
Denote $\theta=PB_-P\phi-P\phi$, the question of the existence of the solution is reduced to studying under which condition $\theta=0$.

We have $P\theta=\theta$ and Theorem~\ref{theo:invertibility} states that $\pi_-\theta=0$. This implies that $\theta\in{\mathcal R}(\pi_+)$. 
The operator ${\mathcal A}^{-1}$ is  symmetric definite positive  on $\mathcal R(\pi_+)$, and induces the scalar product of the ${-1/2}$-norm. 
 Lemma~\ref{lemma::properties::Amoins1} states that, if the constants are controlled   we have $({\mathcal A}^{-1}\theta|\theta)_\HH = 0$, and it follows that $\theta =0$. This proves the result when the constants are controlled.

Assume now that $\Gamma_D\cup\Gamma_R=\emptyset$ and  $\int_\Omega h\ne0$, Lemma~\ref{lemma::properties::Amoins1} states that
\begin{equation}\label{eq:dm1bb}
\Vert \theta\Vert_{-1/2}^2=({\mathcal A}^{-1}\theta|\theta)_\HH=\left(\int_\Omega h\right)^{-1}(\theta|\Phi)_\HH^2.
\end{equation}
If $\int_\Omega h<0$, the two terms have opposite signs, hence both are zero. Then $\theta=0$ and this proves the result for the case   $\Gamma_D\cup\Gamma_R=\emptyset$ and $\int_\Omega h<0$. 

Let us assume now that  $\int_\Omega h>0$. 
Since changing the sign of $\lambda$ amounts to study the same problem where $h$ is replaced by $-h$, we deduce from the case $\int_\Omega h<0$ with $\phi=\Phi$ and the relation $P\Phi=\Phi$ that $PB_+\Phi=\Phi$. Since $\Phi \in \mathcal D(A^{1/2})$, it follows from Theorem \ref{theo:invertibility}  that $B_+\Phi\in \mathcal D(A^{1/2})$. Hence there exists a $s^\star\in H$ such that $B_+\Phi=(1,s^\star)$  and we have $\mathcal AB_+\Phi=(c^{-1}h-c^{-1}\div \sigma \nabla s^\star,0)$. This proves 
$$P\mathcal A B_+ \Phi =\mathcal A B_+ \Phi,$$
 and a simple calculation proves that
\begin{equation}(\mathcal A B_+ \Phi|\Phi)_\HH=(\mathcal A B_+ \Phi|B_+\Phi)_\HH=\int_\Omega h \label{eq:norm}\end{equation}
We then compute 
\[(\Phi-PB_-\Phi|\mathcal A B_+\Phi)_\HH=(\Phi-B_-\Phi|P\mathcal A B_+\Phi)_\HH=(\Phi-B_-\Phi|\mathcal A B_+\Phi)_\HH\underbrace{=}_{\bf (1)}(\Phi|\mathcal A B_+\Phi)_\HH\ne 0,\]
where the equality $\bf (1)$ is obtained by remarking that $ \mathcal  A B_+\Phi\in \mathcal R(\pi_+)$ and  $B_-\Phi \in \mathcal R(\pi_-)$ which are orthogonal spaces. We then obtain $\Phi-PB_-\Phi\ne 0$.

It follows from \eqref{eq:dm1bb} that
\begin{eqnarray*}
\Vert \theta\Vert_{-1/2}^2&=&(\int_\Omega h)^{-1}(\theta|\Phi)_\HH^2=(\int_\Omega h)^{-1}(\theta|PB_+\Phi)_\HH^2\\
&=&(\int_\Omega h)^{-1}(P\theta|B_+\Phi)_\HH^2
=(\int_\Omega h)^{-1}(\theta|B_+\Phi)_\HH^2.
\end{eqnarray*}
Using that $\theta$ and $B_+\Phi$ belong to $\mathcal R(\pi_+)$ on which all the eigenvalues of ${\mathcal A}^{-1}$ are  positive, the above equation implies
$$\Vert \theta\Vert^{2}_{-1/2}=(\int_\Omega h)^{-1}(\theta|{\mathcal A}B_+\Phi)^2_{-1/2}.$$
We recall that $\Vert{\mathcal A} B_+\Phi\Vert_{-1/2}^2=(\mathcal AB_+\Phi|B_+\Phi)_\HH=\int_\Omega h,$
and we obtain 
$$\Vert \theta\Vert^{2}_{-1/2}\Vert{\mathcal A} B_+\Phi\Vert_{-1/2}^2=(\theta|{\mathcal A}B_+\Phi)_{-1/2}^2$$ which is an equality case in Cauchy-Schwarz inequality. This implies that
$\theta$ and $\mathcal AB_+\Phi$ are colinear. Hence there exists some constant $t$ such that \[\theta=t\mathcal AB_+\Phi.\] 

Performing the scalar product with $\Phi$ and using the fact that $(\mathcal AB_+\Phi|\Phi)\neq0$ which follows from \eqref{eq:norm}, we conclude that $t=0$ (hence $\theta=0$)  \FG{if and only if} $(\theta|\Phi)=0$, which reads $(\phi|\Phi-PB_-\Phi)_\HH=0$.

{\em Second case: $ \mathcal K(\mathcal A)\neq\{0\}$, which is the balanced case.}

In the balanced case the kernel of $\mathcal A$ is $\RR \phi_0$, where we recall from section \ref{sec:operatorA} that $P\phi_0=\Phi$.
The condition $\pi_+\psi=0$ is equivalent to the existence of $\alpha \in \RR$ such that $\psi=\pi_-\psi +\alpha \phi_0$. 
The condition $P\psi=P\phi$ is thus equivalent to 
\begin{equation}\label{eq:cns}
P\phi=P\pi_-\psi +\alpha \Phi.
\end{equation}

Necessary condition:

After multiplying Equation \eqref{eq:cns} by $B_-$ one obtains:
\[\pi_-\psi = B_-P\phi-\alpha B_-\Phi.\]
Replacing the expression of $\pi_-\psi$ in \eqref{eq:cns} yields the following necessary condition:
\[PB_-P\phi+\alpha \Phi-\alpha PB_-\Phi=P\phi,\]
which reads
\[\alpha(\Phi-PB_-\Phi)=P\phi-PB_-P\phi.\]
It follows from lemma \ref{lem:lem5} that $(\Phi-PB_-\Phi| \Phi)_\HH\ne 0$, then it is necessary that $$\alpha=\frac{(\Phi-PB_-\Phi|\phi)_\HH}{(\Phi-PB_-\Phi|\Phi)_\HH}.$$

$\psi$ is uniquely determined by 
\begin{equation}\label{eq:defpsi}
\psi=B_-P\phi+\frac{(\Phi-PB_-\Phi|\phi)_\HH}{(\Phi-PB_-\Phi|\Phi)_\HH}( \phi_0-B_-\Phi)
\end{equation}

Conversely, if $\psi$ is defined by \eqref{eq:defpsi}, it is clear that $\pi_+\psi=0$. Let $\theta=P\psi-P\psi$: it suffices to prove that $\theta=0$ to ensure that $\psi$ solves the problem.

$$(\theta|\phi_0)_\HH=(\theta|\Phi)_\HH=0$$
by choice of $\alpha$. A simple calculation shows that
$$\pi_-\theta=0.$$
This proves that $\theta\in \mathcal R(\pi_+)$, where $\mathcal A^{-1}$ is a symmetric positive definite  operator. It follows from lemma \ref{lemma::properties::Amoins1} that $(\mathcal A^{-1}\theta|\theta)_\HH=0$ and hence $\theta=0$. This finishes the proof.

\bibliographystyle{plain} 
\bibliography{biblio}
\end{document}